\documentclass[11pt]{article}
\usepackage{amsmath}
\usepackage{amsfonts}
\usepackage{color}
\usepackage[colorlinks=true]{hyperref}
\usepackage{epsfig}
\usepackage{amssymb}
\def\RIP{{\hbox{\rm RIP}}}
\def\RE{{\hbox{\rm RE}}}
\def\Tr{{\mathop{\hbox{\rm Tr}}}}
\def\Opt{{\mathop{\hbox{Opt}}}}
\def\sign{{\mathop{\hbox{\rm sign}}}}

\def\Card{{\hbox{\rm Card}}}

\def\whg{\widehat{\gamma}_s}
\def\RI{{\hbox{\rm RI}}}

\def\cY{{\cal Y}}

\newcommand{\bbr}{{\Bbb{R}}}

\def\w{\omega}

\def\e{\varepsilon}

\def\Ker {{\rm Ker}}

\newcommand{\epr}{\hfill\hbox{\hskip 4pt
                \vrule width 5pt height 6pt depth 1.5pt}\vspace{0.5cm}\par}

\newcommand{\be}{\begin{eqnarray}}
\newcommand{\ee}[1]{\label{#1}\end{eqnarray}}
\newcommand{\nn}{\nonumber \\}
\newcommand{\ese}{\end{eqnarray*}}
\newcommand{\bse}{\begin{eqnarray*}}
\newcommand{\rf}[1]{~(\ref{#1})}

\newtheorem{lemma}{Lemma}
\newtheorem{proposition}{Proposition}
\newtheorem{corollary}{Corollary}
\newtheorem{theorem}{Theorem}
\newtheorem{definition}{Definition}
\def\argmin{{\mathop{\hbox{\rm argmin}\,}}}
\def\Argmin{{\mathop{\hbox{\rm Argmin}\,}}}

\title{On Verifiable Sufficient Conditions for Sparse Signal Recovery via $\ell_1$ Minimization}
\title{On Verifiable Sufficient Conditions for Sparse Signal Recovery via $\ell_1$ Minimization}
\author{ Anatoli Iouditski\thanks{LJK,
Universit\'e J. Fourier, B.P. 53, 38041 Grenoble
Cedex 9, France, {\tt Anatoli.Juditsky@imag.fr}} \and Arkadi Nemirovski\thanks{Georgia Institute
 of Technology, Atlanta, Georgia
30332, USA, {\tt{nemirovs@isye.gatech.edu}}.
Research of
this author was supported by the Office of Naval Research grant \# N000140811104.}}
\begin{document}
\maketitle
\begin{abstract}
 We discuss necessary and sufficient conditions for a sensing matrix to be ``$s$-good'' -- to allow for exact $\ell_1$-recovery of sparse signals with $s$ nonzero entries when no measurement noise is present. Then we express the error bounds for imperfect $\ell_1$-recovery (nonzero measurement noise, nearly $s$-sparse signal, near-optimal solution of the optimization problem yielding the $\ell_1$-recovery) in terms of the characteristics underlying these conditions. Further, we demonstrate (and this is the principal result of the paper) that these characteristics, although difficult to evaluate, lead to verifiable sufficient conditions for exact sparse $\ell_1$-recovery and to efficiently computable upper bounds on those $s$ for which a given sensing matrix is $s$-good.  We establish also instructive links between our approach and the basic concepts of the Compressed Sensing theory, like Restricted Isometry or Restricted Eigenvalue properties.
\end{abstract}

\section{Introduction}\label{sec:intro}
In the existing literature on sparse signal recovery and Compressed Sensing (see [\ref{r3}-\ref{r9},\ref{r12}-\ref{r15}]  and references therein) the emphasis is on assessing sparse signal $w\in \bbr^n$ from an observation $y\in\bbr^k$ (in this context $k\ll n$):
\be
y=Aw+\xi,\;\;\;\|\xi\|\le \e,
\ee{obs}
where $\|\cdot\|$ is a given norm on $\bbr^k$, $\xi$ is the observation error and $\e\ge 0$ is a given upper bound on the error magnitude, measured in the norm $\|\cdot\|$. One of the most popular
(computationally tractable) estimators which is well suited for recovering sparse signals is the {\sl $\ell_1$-recovery} given by
\be
\widehat{w}
\in\argmin_{z}\left\{\|z\|_1: \;\|Az-y\|\leq\e\right\}.
\ee{L1min}
The existing Compressed Sensing theory focuses on
this estimator and since our main motivation comes from the Compressed Sensing, we will also concentrate on this particular recovery. It is worth to mention that other closely related estimation techniques are used in statistical community, the most renown examples are ``Dantzig Selector'' (cf. \cite{6}), provided by
\be
\widehat{w}'\in\argmin_{z}\left\{\|z\|_1: \;\|A^T(Az-y)\|_\infty\leq\e\right\},
\ee{DSel}
and Lasso estimator, see \cite{Tib,BRT2008},
which under sparsity scenario exhibits similar behavior.
\par
The theory offers strong results which state, in
particular, that if $w$ is $s$-sparse (i.e., has at most {$s$}
nonzero entries) and $A$ possesses a certain well-defined property,
 then the $\ell_1$-recovery of $w$ is close to $w$, provided the observation error $\epsilon$ is small.
  For instance, necessary and sufficient conditions of exactness of $\ell_1$-recovery in the case of noiseless observation (when $\e=0$) has been established in \cite{Zhang,Don1,Don2}. Specifically,
 in
\cite{Zhang} it is shown that $w$ is the unique solution of the noiseless $\ell_1$-recovery problem
\be
\min_{z}\left\{\|z\|_1: \;Az=Aw\right\}.
\ee{noiseless}
if and only if the kernel $\Ker A$ of the sensing matrix is {\em strict $s$-balanced}, the latter meaning that
for any set $I\subset \{1,...,n\}$ of cardinality $\leq s$ it holds
\be
\sum_{i\in I} |z_i|<\sum_{i\not\in I} |z_i|\;\;\mbox{for any}\;\;z\in \Ker A
\ee{zhang1}
(what the above condition is sufficient for the $\ell_1$-recovery to be exact in the noiseless case
was stated in \cite{donhuo}).

 Some particularly impressive results make use of the Restricted Isometry property which is as follows: a $k\times n$ matrix $A$ is said to possess the Restricted Isometry ($\RI(\delta,m)$) property with parameters $\delta\in (0,1)$ and $m$, where $m$ is a positive integer, if
\begin{equation}\label{RIP}
\sqrt{1-\delta}\|x\|_2\leq \|Ax\|_2\leq\sqrt{1+\delta}\|x\|_2 \hbox{\ for all $x\in\bbr^n$ with at most $m$ nonzero entries.}
\end{equation}
For instance,  the following result is well known (\cite[Theorem 1.2]{candescr} or \cite[Theorem
4.1]{ICM}):  let $\|\cdot\|$ in \rf{obs} be the Euclidean norm $\|\cdot\|_2$,
and let the sensing matrix $A$ satisfy $\RI(\delta,2s)$-property with $\delta<\sqrt{2}-1$.
Then
\be
\|\widehat{w}-w\|_1&\le& 2(1-\rho)^{-1}[\alpha\e\sqrt{s}+(1+\rho)\|w-w^s\|_1]
\ee{intro}
where $
\alpha={2\sqrt{1+\delta}\over 1-\delta},\;\;\;\rho={\sqrt{2}\delta\over 1-\delta}$ and $w^s$ is obtained from $w$ by zeroing all but the $s$ largest in absolute values entries.
The conclusion is that when $A$ is $\RI(\delta,2s)$ with $\delta<\sqrt{2}-1$,  $\ell_1$-recovery reproduces well signals with small $s$-tails (small $\|w-w^s\|_1$), provided that the observation error is small.
 Even more impressive is the fact that
there are $k\times n$ sensing matrices $A$ which possess, say, the $\RIP(1/4,{ 2s})$-property for ``large'' {$s$} -- as large as $O\left({k/\ln(n/k)}\right)$. For instance, this is the case, with overwhelming probability, for matrices obtained by normalization (dividing columns by their $\|\cdot\|_2$-norms) of random matrices with i.i.d. standard Gaussian or $\pm1$ entries, as well as for normalizations of random submatrices of the Fourier transform or other orthogonal matrices.
\par
On the negative side, random matrices are the only known matrices which possess the $\RI(\delta,{ 2s})$- property for such large values of {$s$}.
For all known
deterministic families of $k\times n$ matrices provably possessing the $\RI(\delta,{ 2s})$-property, one has
${ s}=O(\sqrt{k})$  (see \cite{deV}), which is essentially worse than the bound ${ s}=O(1)\left(k/\ln(n/k)\right)$ promised by the RI-based theory.
Moreover, $\RI$-property itself is ``intractable'' -- the only currently available technique to verify the {$\RI(\delta,m)$} property for a $k\times n$ matrix amounts to test all its $k\times m$ submatrices. In other words, given a large sensing matrix $A$, one can never be sure that it possesses the RI$(\delta,m)$-property with a given $m\gg1$.
\par
Certainly, the RI-property is not the only property of a sensing matrix $A$ which allows to obtain good error bounds for $\ell_1$-recovery of sparse signals. Two related characteristics are the Restricted Eigenvalue assumption introduced in \cite{BRT2008} and the Restricted Correlation assumption of \cite{Bic}, among others. However, they share with the RI-property not only the nice consequences as in \rf{intro}, but also the drawback of being computationally intractable. To summarize our very restricted and sloppy description of the existing results on $\ell_1$-recovery, neither strict $s$-balancedness, nor Restricted Isometry, or Restricted Correlation assumption and the like, do allow to answer affirmatively the question whether for a {\em given sensing matrix} $A$, an accurate $\ell_1$-recovery of sparse signals with a given number $s$ of nonzero entries is possible.
\par
Now, suppose we face the following problem:  given a sensing matrix ${\cal A}$, which we are allowed to modify in certain ways to obtain a new matrix $A$, our objective is, depending on problem's specifications, either the maximal improvement, or the minimal deterioration of the sensing properties of $A$ with respect to sparse $\ell_1$-recovery. As a simple example, one can think, e.g., of a $2$- or $3$-dimensional  $n$-point grid $E$ of possible locations of signal sources 
and an $N$-element grid $R$ of possible locations of sensors. 
A sensor at a given location measures a known linear form of the signals emitted at the nodes of $E$ which depends on location, and the goal is to place a given number  $k<N$ of sensors at the nodes of $R$ in order to be able to recover, via the $\ell_1$-recovery, all $s$-sparse signals. Formally speaking, we are given  an $N\times n$ matrix ${\cal A}$, and our goal is to extract from it a $k\times n$ submatrix $A$ which is {\sl $s$-good} -- such that whenever the true signal $w$ in (\ref{obs}) is $s$-sparse and there is no observation error ($\xi=0$), the $\ell_1$-recovery (\ref{L1min}) recovers $w$ exactly. To the best of our knowledge, the only existing computationally tractable techniques which allow to approach such a synthesis problem are those based on  {\em mutual incoherence}
\be
\mu(A)=\max_{i\neq j}{|A_i^TA_j|\over A_i^TA_i}
\ee{mu1}
of a $k\times n$ sensing matrix $A$ with columns $A_i$ (assumed to be nonzero).
Clearly, the mutual incoherence can be easily computed even for large matrices. Moreover, bounds of the same type as in \rf{intro} can be obtained for matrices with small mutual incoherence: a matrix $A$ with mutual incoherence $\mu(A)$ and columns $A_j$ of unit $\|\cdot\|_2$-norm satisfies $\RI(\delta,m)$ assumption \rf{RIP} with $\delta=(m-1)\mu(A)$. Unfortunately, the latter relation implies that $\mu$ should be very small to certify the possibility of accurate $\ell_1$-recovery of non-trivial sparse signals, so that the estimates of a ``goodness'' of sensing for $\ell_1$-recovery based on mutual incoherence are very conservative.
\par
 The goal of this paper is to provide new computationally tractable sufficient conditions for sparse recovery.
\par
The overview of our main results is as follows.
\begin{enumerate}
\item Let for $x\in \bbr^n$ $$\|x\|_{s,1}=\max_{\Card(I)\le s}\sum_{i\in I}|x_i|$$ stand for the sum of $s$ maximal magnitudes of components of $x$. Set
\bse
\widehat{\gamma}_s(A)=\max_{x}
\left\{\|x\|_{s,1}:\;\|x\|_1\le 1,\,Ax=0\right\}.
\ese
Starting from optimality conditions for the problem \rf{noiseless} of noiseless $\ell_1$-recovery, we show that $A$ is $s$-good if and only if $\widehat{\gamma}_s(A)<1/2$, thus recovering some of the results of  \cite{Zhang}. While $\widehat{\gamma}_s(A)$ is fully responsible for ideal $\ell_1$-recovery of $s$-sparse signals under {\sl ideal} circumstances, when there is no observation error in (\ref{obs}) and (\ref{L1min}) is solved to precise optimality, in order to cope with the case of imperfect $\ell_1$-recovery (nonzero observation error, nearly $s$-sparse true signal,
(\ref{L1min}) is not solved to exact optimality), we embed the characteristic $\widehat{\gamma}_s(A)$ into a single-parametric family of 
characteristics $\widehat{\gamma}_s(A,\beta)$, $0\leq\beta\leq\infty$. Here
 \bse
\widehat{\gamma}_s(A,\beta)=\max_{x}\left\{\|x\|_{s,1}-\beta\|Ax\|:\; \|x\|_1\leq 1\right\}
 \ese
 (note that $\widehat{\gamma}_s(A,\beta)$ is nonincreasing in $\beta$ and is equal to $\widehat{\gamma}_s(A)$ for all large enough values of $\beta$). We then demonstrate (Section \ref{sec:imperfect}) that whenever $\beta<\infty$ is such that $\widehat{\gamma}_s(A,\beta)<1/2$, the error of  imperfect $\ell_1$-recovery $\widehat{\omega}$ admits an explicit upper bound, similar in structure the RI-based bound (\ref{intro}):
 \[
 \|\widehat{\omega}-\omega\|_1\le (1-2\widehat{\gamma}(A,\beta))^{-1}[2\beta(\e)+2\|w-w^s\|_1+\nu]
 \]
where $\e$ is the measurement error and $\nu$ is the inaccuracy in solving (\ref{L1min}).
\item The characteristics $\widehat{\gamma}_s(A,\beta)$ is still difficult to compute. In
Section \ref{sec:approx},
we develop  efficiently computable lower and upper bounds on $\widehat{\gamma}_s(A,\beta)$. In particular, we show that the quantity $\alpha_s(A,\beta)$,
\[
\alpha_s(A,\beta):=\min\limits_{Y=[y_1,...,y_n]\in\bbr^{k\times n}}\left\{\max_{1\leq j\leq n}
\|(I-Y^TA)e_j\|_{s,1}:\,\|y_i\|_*\leq\beta, 1\leq i\leq n\right\}
\]
(here $\|\cdot\|_*$ is the norm conjugate to $\|\cdot\|$) is an upper bound on $\widehat{\gamma}_s(A,s\beta)$.
\par
 This bound provides us with an efficiently verifiable (although perhaps conservative) {\sl sufficient condition} for $s$-goodness of $A$, namely, $\alpha_s(A,\beta)<1/2$. We demonstrate that our verifiable sufficient conditions for $s$-goodness are less restrictive than those based on mutual incoherence.
 On the other hand, the proposed lower bounds on $\whg(A,\beta)$ allow to bound from above the values of $s$ for which $A$ is $s$-good.

 We also study limitations of our sufficient conditions for
$s$-goodness: unfortunately, it turns out that these conditions, as applied to a $k\times n$ matrix $A$, cannot justify its $s$-goodness when $s>2\sqrt{2k}$, unless $A$ is ``nearly square''. While being much worse than the theoretically achievable, for appropriate $A$'s, level $O(k/\ln(n/k))$ of $s$ for which $A$ may be $s$-good, this ``limit of performance'' of our machinery nearly coincides with the best known values of $s$ for which {\em explicitly given individual $s$-good $k\times n$ sensing matrices}  are known.
\item In Section \ref{sec:RIP}, we investigate the implications of the RI property in our context. While these implications do not contribute to the ``constructive'' part of our results (since the RI property is difficult to verify), they certainly contribute to better understanding of our approach and integrating it into the existing Compressed Sensing theory. The most instructive result of this Section is as follows: whenever $A$ is, say, $\RI(1/4,m)$ (so that the $A$ is $s$-good for $s=O(1)m$),
our verifiable sufficient conditions do certify that $A$ is $O(1)\sqrt{m}$-good -- they guarantee ``at least the square root of the true level $s$ of goodness''.
\item Section \ref{sec:ill} presents some very preliminary numerical illustrations of our machinery. These illustrations, in particular, present experimental evidence of how significantly this machinery can outperform the mutual-incoherence-based one -- the only known to us existing computationally tractable way to certify goodness.
\end{enumerate}
When this paper was finished, we become aware of  the preprint
\cite{dAspremont} which contain results closely related to some of
those in our paper.
The authors of \cite{dAspremont} have ``extracted'' from  \cite{Cohen} the sufficient condition $\whg(A)<1/2$ for $s$-goodness of $A$ and proposed
an efficiently computable upper bound on
$\whg( A)$ based on semidefinite relaxation. This bound
is essentially different from our, and it could be interesting to find
out if one of these bounds is ``stronger'' than the other.

\section{Characterizing $s$-goodness}\label{sec:gamma}
\subsection{Characteristics $\gamma_s(\cdot)$ and $\whg(\cdot)$: definition and basic properties}
The ``minimal'' requirement on a sensing matrix $A$ to be suitable for recovering {\sl $s$-sparse} signals (that is, those with at most $s$ nonzero entries) via $\ell_1$-minimization is as follows: whenever the observation $y$ in (\ref{L1min}) is noiseless and comes from an $s$-sparse signal $w$: $y=Aw$, $w$ should be the unique optimal solution of the optimization problem in (\ref{L1min}) where $\epsilon$ is set to 0. This observation motivates the following
\begin{definition}\label{sgood} Let $A$ be a $k\times n$ matrix and $s$ be an integer, $0\leq s\leq n$.  We say that $A$ is $s$-good, if for every $s$-sparse vector $w\in\bbr^n$, $w$ is the unique optimal solution to the optimization problem
\begin{equation}\label{OptPrb}
\min\limits_{x\in\bbr^n}\left\{\|x\|_1:Ax=Aw\right\}.
\end{equation}
\end{definition}
Let $s_*(A)$ be the largest $s$ for which $A$ is $s$-good; this is a well defined integer, since by trivial reasons every matrix is $0$-good. It is immediately seen that $s_*(A)\leq\min[k,n]$ for every $k\times n$ matrix $A$.
\par
 From now on, $\|\cdot\|$ is the norm on $\bbr^k$
  and $\|\cdot\|_*$ is its conjugate norm:
$$
\|y\|_*=\max_v\left\{v^Ty:\|v\|\leq1\right\}.
$$
We are about to introduce two quantities which are ``responsible'' for $s$-goodness.
\begin{definition}\label{gammas} Let $A$ be a $k\times n$ matrix, $\beta\in[0,\infty]$ and $s\leq n$ be a nonnegative integer. We define the quantities $\gamma_s(A,\beta)$, $\widehat{\gamma}_s(A,\beta)$ as follows:
\par\noindent
{\rm (i)} $\gamma_s(A,\beta)$ is the infinum of  $\gamma\geq0$ such that for every vector $z\in\bbr^n$ with $s$ nonzero entries, equal to $\pm1$, there exists a vector $y\in\bbr^k$ such that
\begin{equation}\label{suchthat}
\|y\|_*\leq\beta\ \&\ (A^Ty)_i\left\{\begin{array}{ll}=z_i,&z_i\neq0\\
\in[-\gamma,\gamma],&z_i=0\\
\end{array}\right.;
\end{equation}
If for some $z$ as above there does not exist $y$ with $\|y\|_*\leq\beta$ such that $A^Ty$ coincides with $z$ on the support of $z$, we set $\gamma_s(A,\beta)=\infty$.
\par\noindent
{\rm (ii)} $\widehat{\gamma}_s(A,\beta)$ is the infinum of  $\gamma\geq0$ such that for every vector $z\in\bbr^n$ with $s$ nonzero entries, equal to $\pm1$, there exists a vector $y\in\bbr^k$ such that
\begin{equation}\label{suchthat1}
\|y\|_*\leq\beta\ \&\ \|A^Ty-z\|_\infty\leq\gamma.
\end{equation}
\end{definition}
To save notation, we will skip indicating $\beta$ when $\beta=\infty$, thus writing $\gamma_s(A)$ instead of $\gamma_s(A,\infty)$, and similarly for $\widehat{\gamma}_s$.
\par
Several immediate observations are in order:
\\
{\bf A}. It is easily seen that the set of the values of $\gamma$ participating in (i-ii) are closed, so that when $\gamma_s(A,\beta)<\infty$, then for every vector $z\in\bbr^n$ with $s$ nonzero entries, equal to $\pm1$, there exists $y$ such that
\begin{equation}\label{suchthat2}
\|y\|_*\leq\beta\ \&\ (A^Ty)_i\left\{\begin{array}{ll}=z_i,&z_i\neq0\\
\in[-\gamma_s(A,\beta),\gamma_s(A,\beta)],&z_i=0\\
\end{array}\right.;\\
\end{equation}
Similarly, for every $z$ as above there exists $\widehat{y}$ such that
\begin{equation}\label{suchthat22}
\|\widehat{y}\|_*\leq\beta\ \&\ \|A^T\widehat{y}-z\|_\infty\leq \widehat{\gamma}_s(A,\beta).\\
\end{equation}
\noindent {\bf B}. The quantities $\gamma_s(A,\beta)$ and $\widehat{\gamma}_s(A,\beta)$ are convex nonincreasing functions of $\beta$, $0\leq\beta<\infty$. Moreover, from {\bf A} it follows that for a given $A$, $s$ and all large enough values of $\beta$ one has $\gamma_s(A,\beta)=\gamma_s(A)$ and $\widehat{\gamma}_s(A,\beta)=\widehat{\gamma}_s(A)$.\\
{\bf C}. Taking into account that the set $\{A^Ty:\|y\|_*\leq\beta\}$ is convex, it follows that if $\gamma_s(A,\beta)<\infty$, then the vectors $y$ satisfying (\ref{suchthat2}) exist for every $s$-sparse vector $z$ with $\|z\|_\infty\leq1$, not only for vectors with exactly $s$ nonzero entries equal to $\pm1$. Similarly, vectors $\widehat{y}$ satisfying (\ref{suchthat22}) exist for all $s$-sparse $z$ with $\|z\|_\infty\leq1$. As a byproduct of these observations, we see that $\gamma_s(A,\beta)$ and $\widehat{\gamma}_s(A,\beta)$ are nondecreasing in $s$.
\par
 Our interest in the quantities $\gamma_s(\cdot,\cdot)$ and $\widehat{\gamma}_s(\cdot,\cdot)$ stems from the following
\begin{theorem} \label{the1} Let $A$ be a $k\times n$ matrix and $s\leq n$ be a nonnegative integer.
\par {\rm (i)} $A$ is $s$-good if and only if $\gamma_s(A)<1$.\par
{\rm (ii)} For every $\beta\in[0,\infty]$ one has
\begin{equation}\label{onehas1}
\begin{array}{ll}
(a)&\gamma:=\gamma_s(A,\beta)<1\Rightarrow \widehat{\gamma}_s\left(A,{1\over1+\gamma}\beta\right)={\gamma\over1+\gamma}<1/2;\\
(b)&\widehat{\gamma}:=\widehat{\gamma}_s(A,\beta)<1/2\Rightarrow {\gamma}_s\left(A,{1\over1-\widehat{\gamma}}\beta\right)={\widehat{\gamma}\over1-\widehat{\gamma}}<1.
\end{array}
\end{equation}
\end{theorem}
The proof of Theorem \ref{the1} is given in Appendix \ref{App1:The1}.
 \par
 Theorem \ref{the1} explains the importance of the characteristic $\gamma_s(\cdot)$
 in the context of $\ell_1$-recovery. However, it is technically more convenient to deal
 with the quantity $\whg(\cdot)$.

\subsection{Equivalent representation of $\widehat{\gamma}_s(A)$}\label{lowerbound}
According to Theorem \ref{the1} (ii), the quantities $\gamma_s(\cdot)$ and $\widehat{\gamma}(\cdot)$ are tightly related. In particular, the equivalent characterization of $s$-goodness in terms of $\whg(A)$ reads as follows:
\[
A \mbox{ is $s$-good}\;\;\Leftrightarrow \;\;\whg(A)<1/2.
\]
In the sequel, we shall heavily utilize an equivalent representation $\whg(A,\beta)$ which, as we shall see in Section \ref{sec:approx}, has important algorithmic consequences. The representation is as follows:
\begin{theorem}\label{widehatgamma} Consider the polytope
\[
P_s=\{u\in\bbr^n:\;\|u\|_1\leq s,\;\|u\|_\infty\leq1\}.
\] One has
\be\widehat{\gamma}_s(A,\beta)=\max_{u,x}\left\{u^Tx-\beta\|Ax\|:\; u\in P_s,\;\|x\|_1\leq 1\right\}.
\ee{onehashuit}
In particular,
\be
\widehat{\gamma}_s(A)=\max_{u,x}\left\{u^Tx:\; u\in P_s,\;\|x\|_1\leq 1,\;Ax=0 \right\}.
\ee{onehasSept}
\end{theorem}
{\bf Proof.} By definition, $\widehat{\gamma}_s(A,\beta)$ is the smallest $\gamma$ such that the closed convex set $C_{\gamma,\beta}:=A^TB_\beta+\gamma B$, where  $B_\beta=\{w\in \bbr^k:\,\|w\|_*\le \beta\}$ and $B=\{v\in\bbr^n:\,\|v\|_\infty\leq1\}$, contains all vectors with $s$ nonzero entries, equal to $\pm 1$. This is exactly the same as to say that $C_{\gamma,\beta}$ contains the convex hull of these vectors; the latter is exactly $P_s$.
 Now, $\gamma$ satisfies the inclusion
$P_s\subset C_{\gamma,\beta}$ if and only if for every $x$ the support function of $P_s$ is majorized by that of $C_{\gamma,\beta}$, namely, for every $x$ one has
\be
\max\limits_{u\in P_s}u^Tx\leq \max_{y\in C(\gamma,\beta)}y^Tx&=&\max\limits_{w,v}\left\{x^TA^Tw+\gamma x^Tv:\,\|w\|_*\le \beta,\,\|v\|_\infty\leq1\right\}\nn
&=&\beta \|Ax\|+\gamma\|x\|_1.
\ee{betafin}
with the convention that when $\beta=\infty$, $\beta\|Ax\|$ is $\infty$ or 0 depending on whether $\|Ax\|>0$ or $\|Ax\|=0$.
That is, $P_s\subset C_{\gamma,\beta}$ if and only if
\[
 \max_{u\in P_s} (u^Tx-\beta\|Ax\|)\leq \gamma\|x\|_1.
\]
By homogeneity w.r.t. $x$, it is equivalent to
$$
\max_{u,x}\left\{u^Tx-\beta\|Ax\|: u\in P_s, \|x\|_1\leq1\right\}\leq\gamma.
$$
Thus, $\widehat{\gamma}_s(A)$ is the smallest $\gamma$ for which the concluding inequality takes place, and we arrive at (\ref{onehashuit}), (\ref{onehasSept}). \epr
Recall that for $x\in\bbr^n$, $\|x\|_{s,1}$  is the sum of the $s$ largest magnitudes of entries in $x$, or, equivalently,
$$\|x\|_{s,1}=\max\limits_{u\in P_s} u^Tx.
$$
Combining Theorem \ref{the1}, and Theorem \ref{widehatgamma}, we get the following
\begin{corollary}
\label{forri}
For a matrix $A\in \bbr^{k\times n}$ {one has} $\whg(A)=\max\limits_x\{\|x\|_{s,1}:Ax=0,\|x\|_1\leq 1\}$, $1\le s\le n$. As a result, matrix $A$ is $s$-good if and only if the maximum of $\|\cdot\|_{s,1}$-norms of vectors $x\in \Ker(A)$ with $\|x\|_1=1$ is $<1/2$.
\end{corollary}
{ Note that (\ref{onehashuit}) and (\ref{onehasSept}) can be seen as an equivalent definition
of $\widehat{\gamma}_s(A,\beta)$, and one can easily prove Corollary \ref{forri} without any reference
to Theorem \ref{the1}, and thus without a necessity even to introduce
the characteristic $\gamma_s(A,\beta)$.  However, we believe that from the methodological point of view the result of Theorem \ref{the1} is important, since it reveals the ``true origin'' of the quantities $\gamma_s(\cdot)$ and $\whg(\cdot)$ as the entities coming from the optimality conditions for the problem \rf{OptPrb}}.

\section{Error bounds for imperfect $\ell_1$-recovery via $\widehat{\gamma}$}\label{sec:imperfect}
We have seen that the quantity $\gamma_s(A)$ (or, equivalently, $\whg(A)$) is responsible for $s$-goodness
of a sensing matrix $A$, that is, for the precise $\ell_1$-recovery of an $s$-sparse signal $w$ in the ``ideal case'' when there is no measurement error and the optimization problem (\ref{OptPrb}) is solved to exact optimality.
It appears that the  same quantities control the error of $\ell_1$-recovery in the case when the vector $w\in \bbr^n$ is not $s$-sparse and the problem (\ref{OptPrb}) is not solved to exact optimality. To see this, let $w^s$, $s\le n$, stand for the best, in terms of $\ell_1$-norm, $s$-sparse approximation of $w$. In other words, $w^s$ is the vector obtained from $w$ by zeroing all coordinates except for the $s$ largest in magnitude.
\begin{proposition}\label{w1control}
Let $A$ be a $k\times n$ matrix, $1\le s\le n$ and let $\whg(A)<1/2$ (or, which is the same, $\gamma_s(A)<1$). Let also $x$ be a $\nu$-optimal approximate solution to the problem  {\rm \rf{OptPrb}}, meaning that
\[
Ax=Aw\;\;\mbox{and $\|x\|_1\le \hbox{\rm $\Opt$}(Aw)+\nu$},
\]
where  $\hbox{\rm $\Opt$}(Aw)$ is the optimal value of \rf{OptPrb}. Then
\[
\|x-w\|_1\le {\nu+2\|w-w^s\|_1\over 1-2\whg(A)}={1+\gamma_s(A)\over 1-\gamma_s(A)}[\nu+2\|w-w^s\|_1].
\]
\end{proposition}
{\bf Proof.}  Let  $z=x-w$ and let $I$ be the set of indices of $s$ largest elements of $w$ (i.e., the support of $w^s$). Denote  by $x^{(s)}$ ($z^{(s)}$) the vector, obtained from $x$ ($z$) by replacing by zero all coordinates of $x$ ($z$) with the indices outside of $I$. As $Az=0$, by Corollary \ref{forri},
\[\|z^{(s)}\|_{1}\le \|z\|_{s,1}\le \whg(A) \|z\|_1.\]
On the other hand, $w$ is a feasible solution to \rf{OptPrb}, so $\hbox{\rm $\Opt$}(Aw)\le \|w\|_1$, whence
\[
\|w\|_1+\nu\ge \|w+z\|_1=\|w^s+z^{(s)}\|_1 +\|(w-w^s)+(z-z^{(s)})\|_1\ge
\|w^s\|_{1}-\|z^{(s)}\|_1+\|z-z^{(s)}\|_1-\|w-w^s\|_1,
\]
or, equivalently,
$$
\|z-z^{(s)}\|_1\leq \|z^{(s)}\|_1+2\|w-w^s\|_1+\nu.
$$
Thus,
\bse
\|z\|_1&=&\|z^{(s)}\|_1+\|z-z^{(s)}\|_1\le 2\|z^{(s)}\|_1+2\|w-w^s\|_1+\nu\\
&\le &2\whg(A) \|z\|_1+2\|w-w^s\|_1+\nu,
\ese
and, as $\whg(A)<1/2$,
\[
\|z\|_1\le {2\|w-w^s\|_1+\nu\over 1-2\whg(A)}.
\]
\epr
We switch now to the properties of approximate solutions $x$ to the problem
\be
\Opt(y)=\min\limits_{x\in\bbr^n}\left\{\|x\|_1:\|Ax-y\|\le \e\right\}
\ee{inexactp}
where $\epsilon\geq0$ and \[
y=Aw+\xi,\;\;\;\xi\in \bbr^k,
\]
 with $\|\xi\|\le \e$.
We are about to show that
in the ``non-ideal case'', when $w$ is ``nearly $s$-sparse'' and \rf{inexactp} is solved to near-optimality, the error of the $\ell_1$-recovery remains ``under control'' -- it admits an explicit upper bound governed by $\gamma_s(A,\beta)$ with a finite $\beta$. The corresponding result is as follows:
\begin{theorem}\label{the2} Let $A$ be a $k\times n$ matrix, $s\leq n$ be a nonnegative integer, let $\epsilon\geq0$,  and let $\beta\in[0,\infty)$ be such that $\widehat{\gamma}:=\whg(A,\beta)<1/2$. Let also $w\in\bbr^n$, let $y$ in  {\rm (\ref{inexactp})} be such that $\|Aw-y\|\leq\epsilon$, and let  $w^s$  be the vector obtained from $w$ by zeroing all coordinates except for the $s$ largest in magnitude. Assume, further, that $x$ is a $(\upsilon,\nu)$-optimal solution to {\rm (\ref{inexactp})}, meaning that
\be
\|Ax-y\|\leq\upsilon\;\; \mbox{and}\;\; \|x\|_1\leq \hbox{\rm $\Opt$}(y)+\nu.
\ee{mv33}
 Then
\be
\|x-w\|_1\leq(1-2\widehat{\gamma})^{-1}[2\beta(\upsilon+\e)+2\|w-w^s\|_1+\nu].
\ee{then332}
\end{theorem}
{\bf Proof.} Since $\|Aw-y\|\leq\epsilon$, $w$ is a feasible solution to (\ref{inexactp}) and therefore $\Opt(y)\leq\|w\|_1$, whence, by (\ref{mv33}),
\begin{equation}\label{whence11}
\|x\|_1\leq\nu+\|w\|_1.
\end{equation}
Let $I$ be the set of indices of
entries in $w^s$. As in the proof of Proposition \ref{w1control} we denote
by $z=x-w$ the error of the recovery, and by $x^{(s)}$ ($z^{(s)}$) the vector obtained from $x$ ($z$) by replacing by zero all coordinates of $x$ ($z$) with the indices outside of $I$.
By \rf{onehashuit} we have
\be
\|z^{(s)}\|_1\le \|z\|_{s,1}\le \beta\|Az\|+\widehat{\gamma}\|z\|_1\le \beta(\upsilon+\e)+\widehat{\gamma}\|z\|_1.
\ee{useful}
On the other hand, exactly in the same way as in the proof of Proposition \ref{w1control} we conclude that
\[
\|z\|_1\leq 2\|z^{(s)}\|_1+2\|w-w^s\|_1+\nu,
\]
which combines with \rf{useful} to imply that
\[
\|z\|_1\le 2\beta(\upsilon+\e)+2\widehat{\gamma}\|z\|_1+2\|w-w^s\|_1+\nu.
\]
Since $\widehat{\gamma}=\whg(A,\beta)<1/2$, this results in
\bse
\|z\|_1\le (1-2\widehat{\gamma})^{-1}[2\beta(\upsilon+\e)+2\|w-w^s\|_1+\nu],
\ese
which is \rf{then332}.\epr
The bound \rf{then332} can be easily rewritten in terms of
$\gamma_s\left(A,{\beta\over 1-\widehat{\gamma}}\right)={\widehat{\gamma}\over1-\widehat{\gamma}}<1$ instead of $\widehat{\gamma}=\whg(A,\beta)$.
\par
The error bound  (\ref{then332}) for imperfect $\ell_1$-recovery, while being in some respects weaker than the RI-based bound (\ref{intro}), is of the same structure as the latter bound: assuming $\beta<\infty$ and $\whg(A,\beta)<1/2$ (or, equivalently, $\gamma_s(A,2\beta)<1$), the error of imperfect $\ell_1$-recovery can be bounded in terms of $\whg(A,\beta)$, $\beta$, measurement error $\epsilon$, ``$s$-tail'' $\|w^s-w\|_1$ of the signal to be recovered and the inaccuracy $(\upsilon,\nu)$ to which the estimate solves the program (\ref{inexactp}). The only flaw in this interpretation is that we need $\whg(A,\beta)<1/2$, while the ``true'' necessary and sufficient condition for $s$-goodness is $\whg(A)<1/2$. As we know, $\whg(A,\beta)=\whg(A)$ for all finite ``large enough'' values of $\beta$, but we do not want the ``large enough'' values of $\beta$ to be really large, since the larger $\beta$ is, the worse is the error bound (\ref{then332}). Thus, we arrive at the question ``what is large enough'' in our context. Here are two simple results in this direction.
\begin{proposition}\label{verynewprop} Let $A$ be a $k\times n$ sensing matrix of rank $k$.\\
{\rm (i)} Let $\|\cdot\|=\|\cdot\|_2$. Then  for every nonsingular $k\times k$ submatrix $\bar{A}$ of $A$ and every $s\leq k$ one has
\begin{equation}\label{nonsingular}
\beta\geq\bar{\beta}=\sigma^{-1}(\bar{A})\sqrt{k},\gamma_s(A)<1\Rightarrow \gamma_s(A,\beta)=\gamma_s(A),
\end{equation}
where $\sigma(\bar{A})$ is the minimal singular value of $\bar{A}$.\\
{\rm (ii)} Let $\|\cdot\|=\|\cdot\|_1$, and let for certain $\rho>0$ the image of the unit $\|\cdot\|_1$-ball in $\bbr^n$ under the mapping $x\mapsto Ax$ contain the ball $B=\{u\in \bbr^k:\,\|u\|_1\le \rho\}$. Then for every $s\leq k$
\begin{equation}\label{ell1}
\beta\geq\bar{\beta}={1\over\rho},\gamma_s(A)<1\Rightarrow \gamma_s(A,\beta)=\gamma_s(A)
\end{equation}
\end{proposition}
{\bf Proof.} Given $s$, let $\gamma=\gamma_s(A)<1$, so that for every vector $z\in\bbr^n$ with $s$ nonzero entries, equal to $\pm1$, there exists $y\in\bbr^k$ such that $(A^Ty)_i=\sign(x_i)$ when $x_i\neq0$ and $|(A^Ty)_i|\leq\gamma$ otherwise. All we need is to prove that in the situations of (i) and (ii) we have $\|y\|_*\leq \bar{\beta}$.\\
In the case of (i) we clearly have $\|\bar{A}^Ty\|_2\leq\sqrt{k}$, whence $\|y\|_*=\|y\|_2\leq \sigma^{-1}(\bar{A})\|\bar{A}^Ty\|_2\leq \sigma^{-1}(\bar{A})\sqrt{k}=\bar{\beta}$, as claimed. In the case of (ii) we have $\|A^Ty\|_\infty\leq1$, whence
$$
\begin{array}{l}
1\geq \max_v\left\{v^TA^Ty:v\in\bbr^n,\|v\|_1\leq1\right\}=\max_u\left\{y^Tu:u=Av,\|v\|_1\leq1\right\}\\
\underbrace{\geq}_{(*)}
\max_u\left\{u^Ty:u\in\bbr^k,\|u\|_1\leq\rho\right\}=\rho\|y\|_\infty=\rho\|y\|_*,\\
\end{array}
$$
where $(*)$ is due to the inclusion $\{u\in\bbr^k:\|u\|_1\leq\rho\}\subset A\left\{v\in\bbr^n:\|v\|_1\leq1\right\}$ assumed in (ii). The resulting inequality implies that $\|y\|_*\leq1/\rho$, as claimed. \epr
\section{Efficient bounding of
$\gamma_s(\cdot)$}
\label{sec:approx}

In the previous section we have seen that the properties of a matrix $A$ relative to $\ell_1$-recovery  are governed by the quantities $\whg(A,\beta)$ -- the less they are, the better. While these quantities is difficult to compute, we are about to demonstrate -- and this is the primary goal of our paper -- that $\whg(A,\beta)$ admits efficiently computable ``nontrivial'' upper and lower bounds.

\subsection{Efficient lower bounding of $\widehat{\gamma}_s(A,\beta)$}\label{sec:lowb}
Recall that $\whg(A,\beta)\ge \whg(A)$ for any $\beta>0$. Thus, in order to provide a lower bound for $\whg(A,\beta)$ it suffices to supply such a bound for $\whg(A)$.
Theorem \ref{widehatgamma} suggests the following scheme
for bounding $\widehat{\gamma}_s(A)$ from below. By (\ref{onehasSept}) we have
$$
\widehat{\gamma}_s(A)=\max_{u\in P_s} f(u),\;f(u)=\max_x\left\{x^Tu: \|x\|_1\leq1,Ax=0\right\}.
$$
Function $f(u)$ clearly is convex and efficiently computable: given $u$ and solving the LP problem
\[x_u\in\hbox{Argmax}_x\left\{u^Tx:\|x\|_1\leq1,Ax=0\right\},
 \]
 we get a linear form $x_u^Tv$ of $v\in P_s$ which underestimates $f(v)$ everywhere and coincides with $f(v)$ when $v=u$. Therefore the easily computable quantity $\max_{v\in P_s}x_u^Tv$ is a lower bound on $\widehat{\gamma}_s(A)$. We now can use the standard {\sl sequential convex approximation} scheme for maximizing the convex function $f(\cdot)$ over $P_s$. Specifically, we run the recurrence
$$
u_{t+1}\in\hbox{Argmax}_{v\in P_s} x_{u_t}^Tv,\;\;\;u_1\in P_s,
$$
thus obtaining a nondecreasing sequence of lower bounds $f(u_t)=x_{u_t}^Tu_t$ on $\widehat{\gamma}_s(A)$. We can terminate this process when the improvement in bounds falls below a given threshold, and can make several runs starting from randomly chosen points $u_1$.
\subsection{Efficient upper bounding of $\widehat{\gamma}_s(A,\beta)$.}
We have seen that the representation (\ref{onehasSept}) suggests a computationally tractable scheme for bounding  $\widehat{\gamma}_s(A)$ from below. In fact, the same representation allows for a tractable way to bound $\widehat{\gamma}_s(A)$ from above, which is as follows.
Whatever be a $k\times n$ matrix $Y$, we clearly have
$$
\max\limits_{u,x}\left\{u^Tx: \|x\|_1\leq1,Ax=0,u\in P_s\right\}=\max\limits_{u,x}\left\{u^T(x-Y^TAx):
\|x\|_1\leq1,Ax=0,u\in P_s\right\},
$$
whence also
$$
\max\limits_{u,x}\left\{u^Tx: \|x\|_1\leq1,Ax=0,u\in P_s\right\}\leq
\max\limits_{u,x}\left\{u^T(x-Y^TAx):
\|x\|_1\leq1,u\in P_s\right\}.
$$
The right hand side in this relation is easily computable, since the objective in the right hand side problem is linear in $x$, and the domain of $x$ in this problem is the convex hull of just $2n$ points $\pm e_i$, $1\leq i\leq n$, where $e_i$ are the basic orths:
\bse
 \max\limits_{u,x}\left\{u^T(x-Y^TAx):\;\|x\|_1\leq1,u\in P_s\right\}
&=&\max\limits_{{u,i,\,\atop1\leq i\leq n}} \left\{|u^T(I-Y^TA)e_i|:\;u\in P_s\right\}\\
&=&\max\limits_{1\leq i\leq n}
\max\limits_{u\in P_s}|u^T(I-Y^TA)e_i|=\max\limits_i\|(I-Y^TA)e_i\|_{s,1}.\\
\ese
 Thus, for all $Y\in\bbr^{k\times n}$,
\bse \widehat{\gamma}_s(A)&=&\max\limits_{u,x}\left\{u^Tx: \|x\|_1\leq1,Ax=0,u\in P_s\right\}\\
&\leq&
f_{A,s}(Y):=\max\limits_{u\in P_s}u^T[(I-Y^TA)e_i]=\max\limits_i\|(I-Y^TA)e_i\|_{s,1},
\ese
so that when setting $\alpha_s(A,\infty):=\min\limits_Y f_{A,s}(Y)$, we get
\bse
\widehat{\gamma}_s(A)\leq\alpha_s(A,\infty).
\ese
Since $f_{A,s}(Y)$ is an easy-to-compute convex function of $Y$, the quantity $\alpha_s(A,\infty)$ also is easy to compute (in fact, this is the optimal value in an explicit LP program with sizes polynomial in $k,n$).
\par
This approach can be easily modified to provide an upper bound for $\widehat{\gamma}_s(A,\beta)$. Namely,
given a $k\times n$ matrix $A$ and
$s\leq k$, $\beta\in[0,\infty]$, let us set
\begin{equation}\label{program}
\alpha_s(A,\beta)=\min\limits_{Y=[y_1,...,y_n]\in\bbr^{k\times n}}\left\{\max_{1\leq j\leq n}
\|(I-Y^TA)e_j\|_{s,1}:\|y_i\|_*\leq\beta, 1\leq i\leq n\right\}.
\end{equation}
As with $\gamma_s$, $\widehat{\gamma}_s$ we shorten the notation $\alpha_s(A,\infty)$ to $\alpha_s(A)$.
\par
 It is easily seen that the optimization problem is (\ref{program}) is solvable, and that $\alpha_s(A,\beta)$ is nondecreasing in $s$, convex and nonincreasing in $\beta$, and is such that $\alpha_s(A,\beta)=\alpha_s(A)$ for all large enough values of $\beta$ (cf. similar properties of $\widehat{\gamma}_s(A,\beta)$). The central  observation in our context is  that {\sl  $\alpha_s(A,\beta)$ is an efficiently computable upper bound on $\widehat{\gamma}_s(A,s\beta)$}, provided that the norm $\|\cdot\|$ is efficiently computable. Indeed, the efficient computability of $\alpha_s(A,\beta)$ stems from the fact that this is the optimal value in an explicit convex optimization problem with efficiently computable objective and constraints. The fact that $\alpha_s$ is an upper bound on $\widehat{\gamma}_s$ is stated by the following
\begin{theorem}\label{aprop2} One has $\widehat{\gamma}_s(A,s\beta)\leq \alpha_s(A,\beta)$.
\end{theorem}
{\bf Proof.} Let $I$ be a subset of $\{1,...,n\}$ of cardinality $\leq s$, $z\in\bbr^n$ be a $s$-sparse vector with nonzero entries equal to $\pm 1$, and let $I$ be the support of $z$. Let $Y=[y_1,...,y_n]$ be such that $\|y_i\|_*\leq\beta$ and
the columns in $\Delta=I-Y^TA$ are of the $\|\cdot\|_{s,1}$-norm not exceeding $\alpha_s(A,\beta)$. Setting $y=Yz$, we have $\|y\|_*\leq \beta\|z\|_1\leq\beta s$ due to $\|y_j\|_*\leq \beta$ for all $j$. Besides this,
$$
\|z-A^Ty\|_\infty=\|(I-A^TY)z\|_\infty=\|\Delta^Tz\|_\infty\leq \alpha_s(A,\beta),
$$
since the $\|\cdot\|_{s,1}$-norms of rows in $\Delta^T$ do not
exceed $\alpha_s(A,\beta)$ and $z$ is an $s$-sparse vector with
nonzero entries $\pm1$. We conclude that
$\widehat{\gamma}_s(A,s\beta)\leq \alpha_s(A,\beta)$, as claimed. \epr
Some comments are in order.
\paragraph{A.}  By the same reasons as in the previous section, it is important to know how large should be $\beta$ in order to have $\alpha_s(A,\beta)=\alpha_s(A)$. Possible answers are as follows. Let $A$ be a $k\times n$ matrix of rank $k$. Then
\\
{\rm (i)} Let $\|\cdot\|=\|\cdot\|_2$. Then  for every nonsingular $k\times k$ submatrix $\bar{A}$ of $A$ and every $s\leq k$ one has
\begin{equation}\label{nonsingularA}
\beta\geq\bar{\beta}={3\over 2}\sigma^{-1}(\bar{A})\sqrt{k},\alpha_s(A)<1/2\Rightarrow \alpha_s(A,\beta)=\alpha_s(A),
\end{equation}
where $\sigma(\bar{A})$ is the minimal singular value of $\bar{A}$.\\
{\rm (ii)} Let $\|\cdot\|=\|\cdot\|_1$, and let for certain $\rho>0$ the image of the unit $\|\cdot\|_1$-ball in $\bbr^n$ under the mapping $x\mapsto Ax$ contain the centered at the origin $\|\cdot\|_1$-ball of radius $\rho$ in $\bbr^k$. Then for every $s\leq k$
\begin{equation}\label{ell1A}
\beta\geq\bar{\beta}={3\over2\rho},\alpha_s(A)<1/2\Rightarrow \alpha_s(A,\beta)=\alpha_s(A)
\end{equation}
The proof is completely similar to the one of Proposition \ref{verynewprop}.
\par
Note that the above bounds on $\beta$ ``large enough to ensure $\alpha_s(A,\beta)=\alpha_s(A)$'', same as their counterparts in Proposition \ref{verynewprop}, whatever conservative they might be, are ``constructive'': to use
the bound (\ref{nonsingularA}), it suffices to find a (whatever) nonsingular $k\times k$ submatrix of $A$ and to compute its minimal singular value. To use (\ref{ell1A}), it suffices to solve $k$ LP programs
$$
\rho_i=\max_{x,\rho}\left\{\rho: \|x\|_1\leq 1, (Ax)_j=\rho\delta_i^j,1\leq j\leq k\right\},\,i=1,...,k
$$
($\delta_i^j$ are the Kronecker symbols) and to set $\rho=\min_i\rho_i$.

\paragraph{B.} Whenever $s$, $t$ are positive integers, we clearly have $\|z\|_{st,1}\leq s\|z\|_{t,1}$, whence
\be
\alpha_{s}(A,\beta)\leq s\alpha_1(A,\beta).
 \ee{as}

 Thus, we can replace in Theorem \ref{aprop2} the quantity $\alpha_s(A,\beta)$ with $s\alpha_1(A,\beta)$.
 Further, we have
 $
 \alpha_1(A,\beta)=\max_i\alpha^i,
 $
 where
\be
\alpha^i:=\min_{y_i}\left\{\|e_i-A^Ty_i\|_\infty:\,\|y_i\|_*\leq\beta\right\},\,i=1,...,n.
\ee{scp}
On the other hand, we have
\bse
 \alpha^i&=&\min_{y}\max_{j}\left\{[e_i-A^Ty]_j:\,\|y_i\|_*\leq\beta\right\}=
 \min_{y}\max_{x}\left\{(e_i-A^Ty)^Tx:\,\|y_i\|_*\leq\beta,\,\|x\|_1\le 1\right\}\\
 &=&\max_{x}\min_{y}\left\{e_i^Tx-y^TAx:\,\|y_i\|_*\leq\beta,\,\|x\|_1\le 1\right\}=
 \max_x \left\{e_i^Tx-\beta \|Ax\|:\,\|x\|_1\le 1\right\}\le \widehat{\gamma}_1(A,\beta),
\ese
and by Theorem \ref{aprop2} we conclude that
\begin{equation}\label{eq444}
\alpha_1(A,{\beta})=\widehat{\gamma}_1(A,\beta),
\end{equation}
i.e. the relaxation for $\widehat{\gamma}_1(A,\beta)$ is exact.
 As a compensation for increased conservatism of the bound \rf{as}, note that while both $\alpha_s$ and $\alpha_1$ are efficiently computable, the second quantity is computationally ``much cheaper''. Indeed, computing $\alpha_1(A,\beta)$ reduces to solving
$n$ convex programs \rf{scp}
of design dimension $k$ each. 
In contrast to this, solving {\rm (\ref{program})} with $s\geq2$ seemingly cannot be decomposed in the aforementioned manner, while ``as it is'' {\rm (\ref{program})} is a convex program of the design dimension $kn$. Unless $k,n$ are small, solving a single optimization program of design dimension $kn$ usually is much more demanding computationally than solving $n$ programs of similar structure with design dimension $k$ each.
\subsection{Relation to the mutual incoherence condition}
\label{sec:MIP}
Remarks in {\bf B}  point at some simple, although instructive conclusions. Let $A$ be a $k\times n$ matrix with nonzero  columns $A_j$, $1\leq j\leq n$, and let $\mu(A)$ be its mutual incoherence, as defined in \rf{mu1}.\footnote{$^)$Note that the ``Euclidean origin'' of the mutual incoherence is not essential in the following derivation. We could start with an arbitrary say, differentiable outside of the origin, norm $p(\cdot)$ on $\bbr^k$, define $\mu(A)$ as $\max\limits_{i\neq j}|A_j^Tp'(A_i)|/p(A_i)$ and define $\beta(A)$ as $\max\limits_{i} \|p'(A_i)/p(A_i)\|_*$, arriving at the same results.}$^)$
\begin{proposition}
\label{extraprop}
For $\beta(A)=\max\limits_{1\leq j\leq n}{\|A_j\|_*\over \|A_j\|_2^2}$ we have
\begin{equation}\label{MIP}
\alpha_1\left(A,{\beta(A)\over1+ \mu(A)}\right)\leq {\mu(A)\over 1+ \mu(A)}.
\end{equation}
In particular, when $\mu(A)<1$ and $1\le s<{1+\mu(A)\over 2\mu(A)}$, one has
\begin{equation}\label{translates}
\gamma_s(A,2s\beta(A))\leq \gamma_s\left(A,{s\beta(A)(1+\mu(A))\over 1-(s-1)\mu(A)}\right)\leq {s\mu(A)\over1-(s-1)\mu(A)}<1.
\end{equation}
\end{proposition}
{\bf Proof.}
Indeed, with $Y_*=[A_1/\|A_1\|_2^2,...,A_n/\|A_n\|_2^2]$, the diagonal entries in $Y_*^TA$ are equal to 1, while the off-diagonal entries are in absolute values $\leq\mu(A)$; besides this, the $\|\cdot\|_*$-norms of the columns of $Y_*$ do not exceed $\beta(A)$. Consequently, for $Y_+={1\over1+\mu(A)}Y_*$, the absolute values of all entries in $I-Y_+^TA$ are $\leq {\mu(A)\over 1+\mu(A)}$, while the $\|\cdot\|_*$-norms of  columns of $Y_+$ do not exceed ${\beta(A)\over 1+\mu(A)}$. We see that the right hand side in the relation
$$
\alpha_1\left(A,{\beta(A)\over 1+\mu(A)}\right)=\min\limits_{Y=[y_1,...,y_n]}\left\{\max\limits_{i,j}|(I-Y^TA)_{ij}|: \|y_i\|_*\leq{\beta(A)\over 1+\mu(A)}\right\}
$$
does not exceed ${\mu(A)\over 1+ \mu(A)}$, since $Y_+$ is a feasible solution for the optimization program in right hand side. This implies the bound \rf{MIP}.

To show \rf{translates} note that from \rf{MIP} with $\beta={\beta(A)\over 1+\mu(A)}$ and \rf{as} we have
\[
\alpha_s(A,\beta)\leq s\alpha_1(A,\beta)\leq {s\mu(A)\over 1+\mu(A)},
 \]
 and it remains to invoke Theorem \ref{aprop2} and Theorem \ref{the1} (ii). \epr

 Observe that taken along with Theorem \ref{the2}, bound (\ref{translates})
 recovers some of the results from    \cite{Donohoetal}.

Proposition \ref{extraprop} implies that computing $\alpha_s(\cdot,\cdot)$ allows to infer that a given $k\times n$ matrix $A$ is $s$-good, for ``reasonably large'' values of $s$. Indeed, take a realization of a random $k\times n$ matrix with independent entries taking values $\pm1/\sqrt{k}$ with probabilities $1/2$. For such a matrix $A$, with an appropriate absolute constant $O(1)$ one clearly has $\mu(A)\leq O(1)\sqrt{\ln(n)/k}$ with probability $\geq1/2$, meaning that $\gamma_s(A,2s\beta(A))\leq 1/2$ for $s\leq O(1)\sqrt{k/\ln(n)}$. Note that such verifiable sufficient conditions for  $s$-goodness based on mutual incoherence are certainly not new, see \cite{Donohoetal}.
We use them here to show that our machinery does allow {\sl sometimes} to justify $s$-goodness for ``nontrivial'' values of $s$, like $O(\sqrt{k/\ln(n)})$.

  \subsection{Application to weighted $\ell_1$-recovery} Note that $\ell_1$-recovery ``as it is'' makes sense only when $A$ is properly normalized, so that, speaking informally, $Ax$ is ``affected equally'' by all entries in $x$. In a general case, one could prefer to use a ``weighted'' $\ell_1$-recovery
\begin{equation}\label{OptPrbScaled}
\tilde{x}_{\Lambda,\epsilon}(y)\in\Argmin_{z\in\bbr^n}\left\{\|\Lambda z\|_1: \|Az-y\|\leq\epsilon\right\},
\end{equation}
where $\Lambda$ is a diagonal matrix with positive diagonal entries $\lambda_i$, $1\leq i\leq n$, which, without loss of generality, we always assume to be $\leq1$. By change of variables $x=\Lambda^{-1}\xi$, investigating $\Lambda$-weighted $\ell_1$-recovery reduces to investigating the standard recovery with the matrix $A\Lambda^{-1}$ in the role of $A$, followed by simple ``translation'' of the results into the language of the original variables. For example, the ``weighted'' version of  our basic Theorem \ref{the2} reads as follows:
\begin{theorem}\label{the22} Let $A$ be a $k\times n$ matrix, $\Lambda$ be a $n\times n$
diagonal matrix with positive entries, $s\leq n$ be a nonnegative integer, and let
$\beta\in[0,\infty)$ be such that $\widehat{\gamma}:=\whg(A\Lambda^{-1},\beta)<1/2$. Let also $w\in\bbr^n$, $\omega=\Lambda w$, and let $\omega^s$ be the vector obtained from $\omega$ by zeroing all coordinates except for the $s$ largest in magnitude. Assume, further, that $y$ in {\rm (\ref{OptPrbScaled})} satisfies the relation $\|Aw-y\|\leq\epsilon$, and that $x$ is a $(\upsilon,\nu)$-optimal solution to {\rm (\ref{OptPrbScaled})}, meaning that
\bse
\|Ax-y\|\leq\upsilon\ \mbox{ and }\ \|\Lambda x\|_1\leq \nu+\hbox{\rm $\Opt$}(y),
\ese
where {\rm$\Opt(y)$} is the optimal value of {\rm (\ref{OptPrbScaled})}. Then
\begin{equation}\label{then3322}
\|\Lambda(x-w)\|_1\leq(1-2\widehat{\gamma})^{-1}[2\beta(\upsilon+\e)+2\|\omega-\omega^s\|_1+\nu].
\end{equation}
\end{theorem}
The issue we want to address here is how to choose the scaling matrix $\Lambda$. When our goal is to recover well signals with as much nonzero entries as possible, we would prefer to make $\whg(A\Lambda^{-1})<1/2$ for as large $s$ as possible (see Theorem \ref{the1}), imposing a reasonable lower bound on the diagonal entries in $\Lambda$ (the latter allows to keep the left hand side in (\ref{then3322}) meaningful in terms of the original variables). The difficulty is that $\whg(A\Lambda^{-1},\beta)$ is hard to compute, not speaking about minimizing it in $\Lambda$. However, we can  minimize in $\Lambda$ the efficiently computable quantity $\alpha_s(A\Lambda^{-1},\bar{\beta})$, $\bar{\beta}=\beta/s$, which is an upper bound on $\widehat{\gamma}_s(A\Lambda^{-1},\beta)$. Indeed, let
\[
\cY=\{Y=[y_1,...,y_n]:\|y_i\|_*\leq\bar{\beta},\;1\le i\le n\}.
\]
 Denoting by $A_i$ the columns of $A$, we have
\bse
\alpha_s(A\Lambda^{-1},\bar{\beta})&=&\min\limits_{Y\in\cY}\left\{\max\limits_{1\leq i\leq n}\|e_i-Y^TA_i\lambda_i^{-1}\|_{s,1}\right\}\\
&=&\min\limits_{Y\in\cY,\,\alpha}\left\{\alpha: \|e_i-Y^TA_i\lambda_i^{-1}\|_{s,1}\leq\alpha,\,1\leq i\leq n\right\}\\
&=&\min\limits_{Y\in\cY,\,\alpha}\left\{\alpha:\|\lambda_ie_i-Y^TA_i\|_{s,1}\leq\alpha\lambda_i,\,1\leq i\leq n\right\},
\ese
so that the problem of minimizing $\alpha_s(A\Lambda^{-1},\bar{\beta})$ in $\Lambda$ under the restriction $0<{\ell}\leq \lambda_i\leq1$ on the diagonal entries of $\Lambda$    reads
\begin{equation}\label{reads44}
\min\limits_{\{\lambda_i\},\,\alpha,\,Y\in\cY}\left\{\alpha:\|\lambda_ie_i-Y^TA_i\|_{s,1}\leq\alpha\lambda_i,\,
{\ell}\leq\lambda_i\leq1, \,1\leq i\leq n\right\}.
\end{equation}
The resulting problem, while not being exactly convex, reduces, by bisection in $\alpha$, to a ``small series'' of explicit convex problems and thus is efficiently solvable. In our context, the situation is even better: basically, all we want is to impose on the would-be $\widehat{\gamma}_s$ an upper bound $\widehat{\gamma}_s(A\Lambda^{-1},\beta)\leq\widehat{\gamma}$
 with a given $\widehat{\gamma}<1/2$, and this reduces to solving a {\sl single} explicit convex feasibility problem
$$
\hbox{find\ } \{\lambda_i\in[{\ell},1]\}_{i=1}^n\hbox{\ and\ } Y\in\cY \hbox{\ such that\ }\|\lambda_ie_i-Y^TA_i\|_{s,1}\leq\widehat{\gamma}\lambda_i,\, 1\leq i\leq n.
$$
\subsection{Limits of performance}
As we have seen in Section \ref{sec:MIP}, the bounding mechanism based on computing $\alpha_s(\cdot,\cdot)$ allows to certify $s$-goodness of an $k\times n$-sensing matrix for $s$ as large as $O(\sqrt{k/\ln(n)})$.  Unfortunately, the $O(\sqrt{k})$-level of values of $s$ is the largest which can be justified via the proposed approach, unless $A$ is ``nearly square''.
\subsection{$\sqrt{k}$-bound}
\begin{proposition}\label{aprop3} For every $k\times n$ matrix $A$ with $n\geq 32k$, every $s$,
$1\leq s\leq n$ and every $\beta\in[0,\infty]$ one has
\begin{equation}\label{eq11}
 \alpha_s(A,\beta)\geq \min\left[{3s\over 4(s+\sqrt{2k})},{1\over 2}\right].
\end{equation}
In particular, in order for $\alpha_s(A,\beta)$ to be $<1/2$ (which, according to Theorems \ref{aprop2} and \ref{the1}, is a verifiable sufficient condition for $s$-goodness of $A$), one should have $s<2\sqrt{2k}$.
\end{proposition}
{\bf Proof.} Let $\alpha:=\alpha_s(A,\beta)$; note that $\alpha\leq1$.
\par
Observe that
\begin{equation}\label{observethat}
\forall v\in\bbr^n: \|v\|_2^2\leq \|v\|_{s,1}^2\max[1,{n\over s^2}].
\end{equation}
Postponing for a while the proof of (\ref{observethat}), let us derive (\ref{eq11}) from this relation. Assume, first, that $s^2\leq n$. Let $Y\in \bbr^{k\times n}$ be such that $\|[I-Y^TA]_j\|_{s,1}\leq \alpha$ for all $j$, where $[B]_j$ is $j$-th column of $B$. Setting $Q=I-Y^TA$, we get a matrix with $\|\cdot\|_{s,1}$-norms of columns $\leq \alpha$. From (\ref{observethat}) it follows that the Frobenius norm of $Q$ satisfies the relation \[
\|Q\|_F^2:=\sum_{i,j}Q_{ij}^2\leq {n^2\alpha^2\over s^2}.\]
 Consequently, \[
 \|Q^T\|_F^2\leq {n^2\alpha^2\over s^2},\]
  whence, setting \[
  H={1\over 2}[Q+Q^T]=I-{1\over 2}[Y^TA+A^TY],\]
   we get \[
   \|H\|_F^2\leq {n^2\alpha^2\over s^2}
    \]as well.  Further, the magnitudes of the diagonal entries in $Q$ (and thus -- in $Q^T$ and in $H$) are at most $\alpha$, whence $\Tr(I-H)\geq n(1-\alpha)$. The matrix $I-H={1\over 2}[Y^TA+A^TY]$ is of the rank at most $2k$ and thus has at most $2k$ nonzero eigenvalues. As we have seen, the sum of these eigenvalues is $\geq n(1-\alpha)$, whence the sum of their squares (i.e., $\|I-H\|_F^2$) is at least ${n^2(1-\alpha)^2\over 2k}$. We have arrived at the relation
$$
{n(1-\alpha)\over \sqrt{2k}}\leq \|I-H\|_F\leq \|I\|_F+\|H\|_F\leq \sqrt{n}+{n\alpha\over s}.
$$
whence
\[
\alpha n\left[{1\over\sqrt{2k}} +{1\over s}\right]\geq {n\over\sqrt{2k}}-\sqrt{n}\geq {3n\over 4\sqrt{2k}}\]
 (the concluding inequality is due to $n\geq 32k$), and (\ref{eq11}) follows. We have derived (\ref{eq11}) from (\ref{observethat}) when $s^2\leq n$; in the case of $s^2>n$, let $s'=\lfloor \sqrt{n}\rfloor$, so that $s'\leq s$. Applying the just outlined reasoning to $s'$ in the role of $s$, we get $\alpha_{s'}(A,\beta)\geq {3s'\over 4(s'+\sqrt{2k})}$, and the latter quantity is $\geq1/2$ due to $n\geq 32k$ and the origin of $s'$. Since $s\geq s'$, we have $\alpha_s(A,\beta)\geq\alpha_{s'}(A,\beta)\geq1/2$, and (\ref{eq11}) holds true.
\par
It remains to prove (\ref{observethat}). W.l.o.g. we can assume that $v_1\geq v_2\geq...\geq v_n\geq0$ and $\|v\|_{s,1}=1$; let us upper bound $\|v\|_2^2$ under these conditions. Setting $v_{s+1}=\lambda$, observe that $0\leq \lambda\leq {1\over s}$ and that for $\lambda$ fixed, we have
$$
\|v\|_2^2\leq \max\limits_{v_1,...,v_s}\left\{\sum_{i=1}^sv_s^2:\sum_{i=1}^sv_i=1,v_i\geq\lambda,\,1\leq i\leq s\right\}+
(n-s)\lambda^2.
$$
The maximum of the right hand side is achieved at an extreme point of the set $\{v\in\bbr^s:\sum_iv_i=1,v_i\geq\lambda\}$, that is, at a point where all but one of $v_i$'s are equal to $\lambda$, and remaining one is $1-(s-1)\lambda$. Thus,
$$\begin{array}{l}
\|v\|_2^2\leq [1-(s-1)\lambda]^2+(s-1)\lambda^2+(n-s)\lambda^2=1-2(s-1)\lambda+(s^2-2s+n)\lambda^2\\
\qquad\qquad\leq\max\limits_{0\leq\lambda\leq 1/s}
[1-2(s-1)\lambda+(s^2-2s+n)\lambda^2].\\
\end{array}
$$
The maximum in the right hand side is achieved at an endpoint of the segment $[0,1/s]$, i.e., is equal to $\max[1,n/s^2]$, as claimed. \epr
\paragraph{Discussion.}
Proposition \ref{aprop3} is a really bad news -- it shows that our verifiable sufficient condition fails to establish $s$-goodness when $s>O(1)\sqrt{k}$, unless $A$ is ``nearly square''. This ``ultimate limit of performance'' is much worse than the actual values of $s$ for which a $k\times n$ matrix $A$ may be $s$-good. Indeed, it is well known, see, e.g. \cite{ICM}, that a random $k\times n$ matrix with i.i.d. Gaussian or $\pm1$ elements is, with close to 1 probability, $s$-good for $s$ as large as $O(1)k/\ln(n/k)$. This is, of course, much larger than the above limit $s\leq O(\sqrt{k})$. Recall, however, that we are interested in {\sl efficiently verifiable}  sufficient condition for $s$-goodness, and efficient verifiability has its price. At this moment we do not know whether the ``price of efficiency'' can be made better than the one for the proposed approach. Note, however, that for all known deterministic provably $s$-good $k\times n$ matrices $s$ is $\leq O(1)\sqrt{k}$, provided $n\gg k$ \cite{deV}.
\section{Restricted isometry property and characterization of $s$-goodness}\label{sec:RIP}
Recall that the RI property (\ref{RIP})
plays the central role in the existing Compressed Sensing results, like the following one: {\sl For properly chosen absolute constants $\delta\in(0,1)$ and integer $\kappa>1$ {\rm (e.g., for $\delta<\sqrt{2}-1$, $\kappa=2$, see \cite[Theorem 1.1]{candescr})}, a matrix possessing  $\RI(\delta,m$) property is $s$-good, provided that $m\geq\kappa s$.}  By Theorem \ref{the1} it follows that with an appropriate $\delta\in(0,1)$, the RI$(\delta,m)$-property of $A$ implies that $\gamma_s(A)<1$, provided $m\geq \kappa s$. Thus, the RI property possesses important implications in terms of the characterization/verifiable sufficient conditions for $s$-goodness as developed above. While these implications do not contribute to the ``constructive'' part of our results (since the RI property is seemingly difficult to verify), they certainly contribute to better understanding of our approach and integrating it into the existing Compressed Sensing theory. In this section, we present the ``explicit forms'' of several of those implications.
\subsection{Bounding $\whg(A)$ for RI sensing matrices}
\begin{proposition}\label{RIprop} Let $s$ be a positive integer,  and let $A$ be a $k\times n$ matrix possessing the {\rm RI$(\delta,2s)$}-property with $0<\delta<\sqrt{2}-1$. Then
\begin{equation}\label{then12}
\whg(A)\leq{{\sqrt{2}\delta\over 1+(\sqrt{2}-1)\delta}}<1/2\;\; \mbox{ and }\;\;
\gamma_s(A)\leq{\sqrt{2}\delta\over 1-\delta}<1.
\end{equation}
\end{proposition}
{\bf Proof.}
Observe that
 by Lemma 2.2 of \cite{candescr}, for any vector $h\in \Ker( A)$ and any index set $I$ of cardinality $\le m/2$ we have under the premise of Proposition:
\[
\sum_{i\in I}|h_i|\le \rho\sum_{i\notin I}|h_i|,\;\;\;\rho={\sqrt{2}\delta(1-\delta)^{-1}}.
\]
This implies that for any $h\in \Ker(A)$ one has $\|h\|_{s,1}\leq\rho(\|h\|_1-\|h\|_{s,1})$, that is, $\|h\|_{s,1}\leq{{\rho\over 1+\rho}}\|h\|_1$. By Corollary \ref{forri} it follows that $\whg(A)\le {{\rho\over 1+\rho}}\,(<1/2)$, and thus $\gamma_s(A)\le \rho\, (<1)$.
\epr
\par
Combining Proposition \ref{RIprop} and Theorem \ref{the1}, we arrive at a sufficient condition for $s$-goodness in terms of RI-property identical to the one in \cite[Theorem 1.1]{candescr}: {\em a matrix $A$ is $s$-good if it possesses the {\rm RI$(\delta,2s)$}-property with $\delta<\sqrt{2}-1$}.
\par
The representation \rf{onehashuit} also allows to bound the value of  $\whg(A,\beta)$ and corresponding $\beta$ in the case when
the Restricted Eigenvalue assumption $\RE(m,\rho,\kappa)$ of \cite{BRT2008} holds true. The exact formulation of the latter assumption is as follows.
Let  $I$ be an arbitrary subset of indices of cardinality $s$; for $x\in \bbr^n$, let  $x^I$ be the vector obtained from $x$ by zeroing all the entries with indices outside of $I$. A sensing matrix $A$ is $\RE(s,\rho,\kappa)$ if
\[
\kappa(s,\rho):=\min_{x,I}\left\{{\|Ax\|_2\over \|x^I\|_{2}}:\;x\in \bbr^n,\,\rho\|x^s\|_{1}\ge\|x-x^s\|_{1};\;{\rm Card}(I)=s
\right\}>0.
\]
Note that the condition $\rho\|x^s\|_{1}\ge \|x-x^s\|_{1}$ is equivalent to $\|x^s\|_{1}\ge {(1+\rho)^{-1}}\|x\|_1$, and
${\|Ax\|_2\over \|x^s\|_2}\ge \kappa$ implies that $\|x^s\|_1\le \kappa^{-1}\sqrt{s}\|Ax\|_2$.
Thus if the $\RE(s,\rho,\kappa)$ assumption holds for $A$, we clearly have for any $x\in \bbr^n$
\bse
\|x\|_{s,1}\le \max\left\{{\sqrt{s}\|Ax\|_2\over \kappa},\;(1+\rho)^{-1}\|x\|_1\right\}.
\ese
In other words,
assumption $\RE(s,\rho,\kappa)$ implies that
\[
\whg\left(A,{\sqrt{s}\over \kappa}\right)\le (1+\rho)^{-1}.
\]
\subsection{``Large enough'' values of $\beta$}
We present here an upper bound on the value of $\beta$ such that $\gamma_s(A,\beta)=\gamma_s(A)$ in the case when the matrix $A$ possesses the RI-property (cf. Proposition \ref{verynewprop}):
\begin{proposition}\label{kappaofA} Let $s$ be a positive integer,  $A$ be a $k\times n$ matrix possessing the {\rm RI$(\delta,2s)$}-property with $0<\delta<\sqrt{2}-1$ and  $s\le n$ and let $\|\cdot\|$ be the $\ell_2$-norm. Then
\be
\widehat{\gamma}_s(A,\beta)\leq
{\sqrt{2}\delta\over 1+(\sqrt{2}-1)\delta}\;\;\mbox{ for all }\; \beta\geq
{\sqrt{(1+\delta)s}\over1+(\sqrt{2}-1)\delta}.
  \ee{onehas678}
\end{proposition}
{\bf Proof.}
The derivations below are rather standard to Compressed Sensing.
Let us prove that
\begin{equation}\label{letusprove}
\forall w\in\bbr^n: \|w\|_{s,1}\leq {\sqrt{(1+\delta)s}\over1+(\sqrt{2}-1)\delta}\|Aw\|_2+
{\sqrt{2}\delta\over 1+(\sqrt{2}-1)\delta}\|w\|_1.
\end{equation}
There is nothing to prove when $w=0$; assuming $w\neq0$, by homogeneity we can assume that $\|w\|_1=1$. Besides this, w.l.o.g. we may assume that $|w_1|\geq|w_2|\geq...\geq |w_n|$.  Let us set $\alpha=\|Aw\|_2$.
 Let us split $w$ into consecutive $s$-element ``blocks'' $w^0,w^1,...,w^q$, so that $w^0$ is obtained from $w$ by zeroing all coordinates except for the first $s$ of them, $w^1$ is obtained from $w$ by zeroing all coordinates except of those with indices $s+1,s+2,...,2s$, and so on, with evident modification for the last block $w^q$.  By construction we have
 \[
 w=\sum_{j=0}^q w^j,\;\;\;\|w^0\|_1\ge \|w^1\|_1\ge ...\ge \|w^q\|_1,\;\; \|w\|_1=\sum_{j=0}^q \|w^j\|_1.
 \]
Further, we have due to the monotonicity of $|w_i|$ and $s$-sparsity of all $w^j$:
\be
j\geq 1\Rightarrow \|w^j\|^2_2\le \|w^j\|_\infty\|w^j\|_1\leq s^{-1}\|w^{j-1}\|_1 \|w^j\|_1\le  s^{-1}\|w^{j-1}\|^2_1 .
\ee{l2l1}
On the other hand, due to the RI-property of $A$ and the fact that $\|Aw\|_2=\alpha$ we have the first inequality in the following chain:

\be
\alpha\sqrt{1+\delta}\|w^0+w^1\|_2&\ge& \|Aw\|_2\|A(w^0+w^1)\|_2\ge (Aw)^TA(w^0+w^1)\nn
&=&(w^0+w^1)^TA^TA(w^0+w^1)+\sum_{j=2}^q (w^0+w^1)^TA^TAw^j\nn
&\ge & (1-\delta)\|w^0+w^1\|_2^2-\sum_{j=2}^q \sqrt{2}\delta \|w^0+w^1\|_2\|w^j\|_2,
\ee{start1}
where we have used the ``classical'' RI-based relation (see \cite{7})
\[
v^TA^TAu\le \sqrt{2}\delta \|v\|_2\|u\|_2
\]
for any two vectors $u,\,v\in \bbr^n$ with disjoint supports and such that $u$ is $s$-sparse and $v$ is $2s$-sparse.  Using \rf{l2l1} we can now
continue \rf{start1} to get
\bse
(1-\delta)\|w^0+w^1\|_2^2&\le& \alpha\sqrt{1+\delta}\|w^0+w^1\|_2+ \sqrt{2}\delta\|w^0+w^1\|_2\;s^{-1/2}\sum_{j=1}^{q-1}\|w^j\|_1\\
&\le & \alpha\sqrt{1+\delta}\|w^0+w^1\|_2+\sqrt{2}s^{-1/2}\delta\|w^0+w^1\|_2\|w-w^0\|_1.
\ese
Since $w^0$ is $s$-sparse,  we conclude that
\[
\|w^0\|_1\le \sqrt{s}\|w^0\|_2\le \sqrt{s}\|w^0+w^1\|_2\le {\alpha\sqrt{(1+\delta)s}
\over 1-\delta}+{\rho} \|w-w^0\|_1= {\alpha\sqrt{(1+\delta)s}
\over 1-\delta}+ \rho(1-\|w^0\|_1)\eqno{[\rho={\sqrt{2}\delta\over1-\delta}]}
\]
(recall that $\|w\|_1=1$). It follows that
\[
\|w^0\|_1\le {\alpha\sqrt{(1+\delta)s}
\over (1+\rho)(1-\delta)}+{\rho\over1+\rho}={\alpha\sqrt{(1+\delta)s}\over1+(\sqrt{2}-1)\delta}+
{\sqrt{2}\delta\over 1+(\sqrt{2}-1)\delta}.
\]
Recalling that $\alpha=\|Aw\|_2$, the concluding inequality is exactly (\ref{letusprove}) in the case of $\|w\|_1=1$. (\ref{letusprove}) is proved.
\par
Invoking \rf{onehashuit}, (\ref{letusprove}) implies that with $\|\cdot\|=\|\cdot\|_2$ and with $\beta\geq{\sqrt{(1+\delta)s}\over1+(\sqrt{2}-1)\delta}$, one has $\widehat{\gamma}_s(A,\beta)\leq
{\sqrt{2}\delta\over 1+(\sqrt{2}-1)\delta}$. \epr
It is worth to note that when using the bounds of Proposition \ref{kappaofA} on $\whg(A,\beta)$ and the corresponding $\beta$ along with Theorem \ref{the2}, we recover the classical bounds on the accuracy of the $\ell_1$-recovery as those given in \cite{ICM, candescr}.

\subsection{Performance of verifiable conditions for $s$-goodness in the case of RI sensing matrices}
It makes sense to ask how conservative is the {\sl verifiable}
sufficient condition for $s$-goodness ``$\alpha_s(A)<1/2$'' as compared to the {\sl difficult-to-verify} RI condition ``if $A$ is $\RI(\delta,m)$, then $A$ is $s$-good for $s\leq O(1)m$''. It turns out that this conservatism is under certain control,
fully compatible with the ``limits of performance'' of our verifiable condition as stated in Proposition \ref{aprop3}. Specifically, we are about to prove that if $A$ is $\RI(\delta,m)$, then $\alpha_s(A)<1/2$ when $s\leq O(1)\sqrt{m}$: our verifiable condition ``guarantees at least square root of what actually takes place''. The precise statement is as follows:
\begin{proposition}\label{aprop} Let a $k\times n$ matrix $A$ possess $\RI(\delta,m)$-property. Then
\begin{equation}\label{thenagain}
\alpha_1(A)\leq
{\sqrt{2}\,\delta\over (1-\delta)\sqrt{m-1}},
\end{equation}
so that
\begin{equation}\label{sothat13}
s<{(1-\delta)\sqrt{m-1}\over 2\sqrt{2}\delta}\;\Rightarrow \;\alpha_s(A)\leq s\alpha_1(A)<1/2.
\end{equation}
\end{proposition}
{\bf Proof.} 1$^0$. We start with the following simple fact (cf. Proposition \ref{RIprop}):
\begin{lemma}\label{lemnew} Let $A$ possess $\RI(\delta,m)$-property. Then \be
\widehat{\gamma}_1(A)\leq
{\sqrt{2}\delta\over (1-\delta)\sqrt{m-1}}.
\ee{g1}
\end{lemma}
{\bf Proof.}
Invoking Theorem \ref{widehatgamma}, all we need to prove is that under the premise of Lemma for every $s$, $1\leq s<m$, and  for every $w\in \Ker(A)$ we have
\[
\|w\|_\infty=\|w\|_{1,1}\leq \widehat{\gamma}\|w\|_1.
\]
To prove this fact we use again the standard machinery related to the RI-property (cf proof of Proposition \ref{kappaofA}): we set $t= \lfloor m/2\rfloor$, assume w.l.o.g. that $\|w\|_1=1$, $|w_1|\geq|w_2|\geq...\geq |w_n|$ and split $w$ into $q$ consecutive ``blocks'' so that the cardinality of the ``blocks'' is $1, \,t-1, t,\,t,\,...$. I.e. the first ``block'' $\w^0\in \bbr^n$ is the vector such that $w^0_1=w_1$ and all other coordinates vanish, $w^1$ is obtained from $w$ by zeroing all coordinates except of those with indices $2,3,...,t$, $w^2$ is obtained from $w$ by zeroing all coordinates except of those with indices $t+1,...,2t$, and so on, with evident modification for the last vector $w^q$.
 Acting as in the proof of Proposition \ref{kappaofA}, and using the relation (see \cite{7})
\[
v^TA^TAu\le \delta \|v\|_2\|u\|_2
\]
for any two $t$-sparse vectors $u,\,v\in \bbr^n$, $t\le m/2$, with disjoint supports, we obtain
 \[
 0=(A(w^0+w^1))^TAw\ge (1-\delta)\|w^0+w^1\|_2^2-t^{-1/2}\delta \|w^0+w^1\|_2
 \]
 whence
 \[
 |w_1|\le \|w^0+w^1\|_2\le {\delta\over (1-\delta)\sqrt{t}},
 \]
what is \rf{g1}. \epr
\par\noindent 2$^0$. Now we are ready to complete the proof of (\ref{thenagain}). We already know that $\alpha_s(A)\leq s\alpha_1(A)$, so all we need is to verify (\ref{thenagain}).  The latter is  readily given by (\ref{eq444}) combined with (\ref{g1}). \epr

\section{Numerical illustration}\label{sec:ill}
We are about to present some very preliminary numerical results for relatively small sensing matrices.
\paragraph{The data.}  In the two series of experiments presented below we deal with sensing matrices of row dimension $n=256$ and $n=1024$.

For $n=256$ we generate three sets of random matrices of column dimension $m=0.1n,\,...,\, 0.9n$: Gaussian matrices, with the i.i.d. normal entries, Fourier matrices, which are $m$ rows of the Fourier basis on $[0,1]$ drawn at random and, finally, Hadamard matrices, which are, again, random $m\times n$ cuts from the $n\times n$ Hadamard matrix.\footnote{Hadamard matrix $H_\ell$ of order $n=2^\ell$ is the orthogonal  matrix with entries $\pm1$ given by the recurrence $H_0=1,$ $H_{\ell+1}=[H_\ell,H_\ell;H_\ell,-H_\ell]$.} Then all matrices are normalized so that their columns have unit $\ell_2$-norm.

For $n=1024$ we provide the result of an experiment with a family of Gaussian matrices of column dimension $m=0.1n,\,...,\, 0.9n$ and with a $992\times 1024$ matrix $A_{\rm conv}$ which is constructed as follows. Let us consider a signal $x$ ``living'' on ${\mathbf{Z}}^2$ and supported on the $32\times 32$ grid $\Gamma=\{(i,j)\in{\mathbf{Z}}^2: \;0\leq i,j\leq 31\}$. We subject such a signal to discrete time convolution with a kernel supported on the set $\{(i,j)\in{\mathbf{Z}}^2: -7\leq i,j\leq 7\}$, and then restrict the result on the $32\times 31$ grid  $\Gamma_+=\{(i,j)\in\Gamma: 1\leq j\leq 31\}$. This way we obtain a linear mapping $x\mapsto A_{\rm  conv}x: \;\bbr^{1024}\to\bbr^{992}$.

\paragraph{The goal} of the experiment is to bound from below and from above the maximal $s$ for which the $m\times n$ matrix $A$ in question is $s$-good (the quantity $s_*(A)$ from Definition \ref{sgood}).
\paragraph{The lower bound} on $s_*(A)$ was obtained via bounding from above, for various $s$, the quantity $\alpha_s(A)$ and invoking Theorem \ref{aprop2} and Theorem \ref{the1} (ii) which, taken together, state that a sufficient condition for $A$ to be $s$-good is $\alpha_s(A)<1/2$.

We provide two lower bounds for $s_*(A)$. The first bound  is obtained using the upper bound $\alpha_s(A)\leq s\alpha_1(A)$ (see Comment {\bf B} in Section \ref{sec:approx}). When the upper bound $s\alpha_1(A)$ for $\alpha_s(A)$ is computed and turns out to be $<1/2$, we know that $A$ is $s$-good, and  our lower bound on $s_*(A)$ is the largest $s$ for which this situation takes place; note that computing this bound reduces to a {\sl single} computation of $\alpha_1(A)$. As explained in Comment {\bf B}, this computation reduces to solving $n$ convex programs of design dimension $m$ each, and these programs are easily convertible to LP's with $(2n+1)\times (m+1)$ constraint matrices. These LP's were solved using the commercial LP solver {\tt mosekopt} \cite{mosek}. Note that in fact computing $\alpha_1(A)$ allows to somehow improve the trivial upper bound $s\alpha_1(A)$ on $\alpha_s(A)$, specifically, as follows. As a result of computing $\alpha_1(A)$, we get the associated matrix $Y$; the largest of $\|\cdot\|_{s,1}$-norms of the columns of $I-Y^TA$ clearly is an upper bound on $\alpha_s(A)$, and this bound is {\sl at worst} $s\alpha_1(A)$.
\par
For ``small'' matrices with the row dimension $n=256$ we also provide the ``improved'' lower bound, obtained using the computation of $\alpha_s(A)$ itself. We act as follows: when the bound $s(\alpha_1)$ is computed, verify if the value $s(\alpha_1)+1$ can be certified lower bound for $s_*(A)$ using the computation of $\alpha_s(A)$. If this bound is certified we proceed with $s(\alpha_1)+2$, and so on. Note that, exactly as it is in the case of $\alpha_1(A)$,  computing $\alpha_s(A)$ allows to improve the lower bound on $s_*(A)$ in the case when $\alpha_s(A)<1/2$. Indeed, as a result of computing $\alpha_s(A)$, we get the associated matrix $Y$; the largest $s$ such that the $\|\cdot\|_{s,1}$-norms of the columns of $I-Y^TA$ is $<1/2$  clearly is a lower bound on $s_*(A)$.
\par
We would like to add here two words about the techniques used to compute the corresponding bound (being of interest  by themselves, these techniques are the subject of a separate paper).  While $\alpha_s(A)$  is efficiently computable via LP (when $\beta=\infty$, the optimization program in  (\ref{program}) is easily convertible into a linear programming one), the sizes of the resulting LP are rather large -- when $A$ is $m\times n$, the LP reformulation of (\ref{program}) has a $(2n^2+n)\times (n(m+n+1)+1)$ constraint matrix. For instance, for $m=230$ and $n=256$, the size of the LP becomes  131,328$\times$127,233, and we preferred to avoid solving this, by no means small, LP program using the interior-point solver available with {\tt mosekopt}. Instead, the LP is reformulated as a saddle-point problem and is solved using an implementation of the non-Euclidean mirror-prox algorithm, described in \cite{mirror}.

\paragraph{The upper bound} on $s_*(A)$ is computed using the lower bound on $\gamma_s(A)$ by the Sequential Convex Approximation algorithm presented in Section \ref{sec:lowb}.
%
%
\paragraph{The results} of our experiments are presented in Tables \ref{tab1} and \ref{tab2}.
The computations we run on an Intel P9500@2.53GHz CPU (the computations were running single-core). We present along with the results the corresponding CPU usage.

 We would like to add the following comment:
our efficiently computable lower bounds on $s_*(A)$ outperform significantly those based on mutual incoherence. Further, these lower and upper bounds ``somehow'' work in the case of the randomly chosen sensing  matrix and work quite well
in the case of the convolution matrix. While the gap between the lower and the upper bound in the case of the random sensing matrix could be better, we can re-iterate at this point our remark that computability has its price.

\begin{table}[h]
\centerline{
\begin{tabular}{||c||c|c|c|c||r|r|r||}
\multicolumn{8}{c}{Gaussian matrix}\\
\hline
&\multicolumn{3}{|c|}{lower bounds on $s_*(A)$}&upper&\multicolumn{3}{|c||}{CPU time (s)}\\
\cline{2-4}\cline{6-8}
$m$&$s[\mu]$&$s[\alpha_1]$&$s[{\alpha}_s]$&bound $\overline{s}$&
\multicolumn{1}{|c|}{$s[\alpha_1]$}&\multicolumn{1}{|c|}{$s[{\alpha}_s]$}&\multicolumn{1}{|c||}{$\overline{s}$}\\
\hline
\hline
25&1&1&1&1&11.0&21.6&3.4\\ \hline
51&1&2&2&4&22.3&24.1&8.8\\ \hline
76 &1&3&3&7&34.2&34.3&23.1\\ \hline
102 &1&3&4&11&50.8&190.7&34.0\\ \hline
128 &1&5&5&15&69.3&75.8& 31.6\\ \hline
153 &1&5&6&19&93.8&557.6&60.7\\ \hline
179 &2&7&8&25&115.4&658.3&103.8\\ \hline
204 &2&9&11&31&141.2&551.5&97.8\\ \hline
230& 2&14&17&41&173.0&561.0& 97.8
\\
\hline
\multicolumn{8}{c}{~}\\
\multicolumn{8}{c}{Random Fourier matrix}\\
\hline
&\multicolumn{3}{|c|}{lower bounds on $s_*(A)$}&upper&\multicolumn{3}{|c||}{CPU time (s)}\\
\cline{2-4}\cline{6-8}
$m$&$s[\mu]$&$s[\alpha_1]$&$s[{\alpha}_s]$&bound $\overline{s}$&\multicolumn{1}{|c|}{$s[\alpha_1]$}&\multicolumn{1}{|c|}{$s[{\alpha}_s]$}&\multicolumn{1}{|c||}{$\overline{s}$}\\
\hline
\hline
24&1&1&1&2&9.3&6.1&1.3\\ \hline
51&1&2&2&4&129.5&14.5&7.2\\ \hline
76 &2&3&3&6&233.1&12.8&16.1\\ \hline
102 &2&4&4&7&213.9&11.2&25.6\\ \hline
128 &2&4&4& 8&270.9&426.5&58.1\\ \hline
152 &3&5&5&10&245.9&2350.7&57.8\\ \hline
178 &3&6&6&14&319.7&161.2&81.5\\ \hline
204 &4&7&7&14&234.0&97.9&75.8\\ \hline
230& 4&9&9& 19&343.2&76.0&51.9
\\
\hline
\multicolumn{8}{c}{~}\\
\multicolumn{8}{c}{Random Hadamard matrix}\\
\hline
&\multicolumn{3}{|c|}{lower bounds on $s_*(A)$}&upper&\multicolumn{3}{|c||}{CPU time (s)}\\
\cline{2-4}\cline{6-8}
$m$&$s[\mu]$&$s[\alpha_1]$&$s[{\alpha}_s]$&bound $\overline{s}$&\multicolumn{1}{|c|}{$s[\alpha_1]$}&\multicolumn{1}{|c|}{$s[{\alpha}_s]$}&\multicolumn{1}{|c||}{$\overline{s}$}\\
\hline
\hline
25&1&1&1&2&10.1&7.4&1.2\\ \hline
51&1&2&2&4&21.6&11.7&3.5\\ \hline
76 &2&3&3&4&34.1&14.2&6.7\\ \hline
102 &3&4&4&11&50.8&23.8&37.7\\ \hline
128 &3&5&5& 7&69.6&48.5&24.1\\ \hline
153 &3&7&7&11&93.8&31.1&84.7\\ \hline
179 &4&9&9&15&112.0&51.0&88.9\\ \hline
204 &5&12&12&15&141.6&51.1&78.6\\ \hline
230& 6&18&18& 28&141.5&55.4&44.1
\\
\hline
\end{tabular}}
\caption{\label{tab1} Efficiently computable bounds on $s_*(A)$, $n=256$.}
{\footnotesize {\bf Lower bound} $s[\mu]$: the bound (\ref{translates}) based on mutual incoherence; $s[\alpha_1]$-bound: the ``improved'' bound based on upper bounding of $\alpha_s(A)$ via the matrix $Y$ obtained when computing $\alpha_1(A)$;  $s[\alpha_s]$: the bound based on computing $\alpha_s(A)$.
{\bf Upper bound} $\overline{s}$: the bound based on successive convex approximation}
\end{table}
%

\begin{table}[h]
\centerline{
\begin{tabular}{||c||c|c|c|c||r|r|r||}
\multicolumn{6}{c}{Gaussian matrix}\\
\hline
&\multicolumn{2}{|c|}{lower bounds on $s_*(A)$}&upper&\multicolumn{2}{|c||}{CPU time (s)}\\
\cline{2-3}\cline{5-6}
$m$&$s[\mu]$&$s[\alpha_1]$&bound $\overline{s}$&\multicolumn{1}{|c|}{$s[\alpha_1]$}&\multicolumn{1}{|c||}{$\overline{s}$}\\
\hline
\hline
102&2&2&8&457.0&400.7\\ \hline
204&2&4&18&1179.0&1722.1\\ \hline
307&2&6&30&2234.6&7585.9\\ \hline
409&3&7&44&3658.6&3421.7\\ \hline
512&3&10&61&5341.7&6304.3\\ \hline
614&3&12&78&7155.7&17616.7\\ \hline
716&3&15&105&9446.1&11670.4\\ \hline
819&4&21&135&12435.1&8373.1\\ \hline
921&4&32&161&13564.2&9838.3\\
\hline
\multicolumn{6}{c}{~}\\
\multicolumn{6}{c}{Convolution matrix}\\
\hline
&\multicolumn{2}{|c|}{lower bounds on $s_*(A)$}&upper&\multicolumn{2}{|c||}{CPU time (s)}\\
\cline{2-3}\cline{5-6}
$m$&$s[\mu]$&$s[\alpha_1]$&bound $\overline{s}$&\multicolumn{1}{|c|}{$s[\alpha_1]$}&\multicolumn{1}{|c||}{$\overline{s}$}\\
\hline
\hline
960&0&5&7&4579.1&271.8
\\
\hline
\end{tabular}}
\caption{\label{tab2} Efficiently computable bounds on $s_*(A)$, $n=1024$.}
{\footnotesize {\bf Lower bound} $s[\mu]$: the bound (\ref{translates}) based on mutual incoherence; $s[\alpha_1]$-bound: the ``improved'' bound based on upper bounding of $\alpha_s(A)$ via the matrix $Y$ obtained when computing $\alpha_1(A)$.
{\bf Upper bound} $\overline{s}$: the bound based on successive convex approximation}
\end{table}

\appendix
\section{Proof of Theorem \ref{the1}} \label{App1:The1}
{\bf Proof.} (i): a) Assume that $A$ is $s$-good, and let us prove that $\gamma_s(A)<1$. Let $I$ be an $s$-element subset of the index set $\{1,...,n\}$ and $\bar{I}$ be its complement, and let $w$ be a vector supported on $I$ and with nonzero $w_i$, $i\in I$. Then $w$ should be the unique solution to the LP problem (\ref{OptPrb}).
From the fact that $w$ is an optimal solution to this problem it follows, by optimality conditions, that for certain $y$ the function
$f_y(x)=\|x\|_1-y^TAx$ attains its minimum over $x\in\bbr^n$ at $x=w$, meaning that $0\in\partial f_y(w)$, that is,
$$
(A^Ty)_i\left\{\begin{array}{ll}=\sign(w_i),&i\in I\\
\in[-1,1],&i\in \bar{I}\\
\end{array}\right.,
$$
so that the LP problem
\begin{equation}\label{one}
\min_{y,\gamma}\left\{\gamma: (A^Ty)_i\left\{\begin{array}{ll}=\sign(w_i),&i\in I\\
\in[-\gamma,\gamma],&i\in\bar{I}\\
\end{array}\right.\right\}
\end{equation}
has optimal value $\leq1$. Let us prove that in fact the optimal value is $<1$. Indeed, assuming that the optimal value is exactly 1, there should exist Lagrange multipliers $\{\mu_i:i\in I\}$ and $\{\nu_i^\pm\geq0:i\in\bar{I}\}$ such that the function
$$
\gamma+\sum_{i\not\in I}\left[\nu_i^+[(A^Ty)_i-\gamma]+\nu_i^-[-(A^Ty)_i-\gamma]\right]-\sum_{i\in I}\mu_i\left[(A^Ty)_i-\sign(w_i)\right]
$$
has unconstrained minimum in $\gamma,y$ equal to 1, meaning that
$$
\begin{array}{ll}
(a)&\sum_{i\in\bar{I}}[\nu_i^++\nu_i^-]=1,\\
(b)&\sum_{i\in I}\mu_i\sign(w_i)=1,\\
(c)&Ad=0, \hbox{ where $d\in\bbr^n$ with\ }d_i=\left\{\begin{array}{ll}-\mu_i,&i\in I\cr
\nu_i^+-\nu_i^-,&i\in\bar{I}.\cr
\end{array}\right. \\
\end{array}
$$
Now consider the vector $x_t=w+td$, where $t>0$. This is a feasible solution to (\ref{OptPrb}) due to $(c)$;
the $\|\cdot\|_1$-norm of this solution is
$$
\sum_{i\in I}|w_i-t\mu_i|+t\sum_{i\in\bar{I}}|\nu_i^+-\nu_i^-|\leq
\sum_{i\in I}|w_i-t\mu_i|+t
$$
where the concluding inequality is given by $(a)$ and the fact that $\nu_i^\pm\geq0$. Since $w_i\neq0$ for $i\in I$,  for small positive $t$ we have
$$
\sum_{i\in I}|w_i-t\mu_i|=\sum_{i\in I}|w_i| -t\sum_{i\in I}\mu_i\sign(w_i)=\sum_{i\in I}|w_i| - t,
$$
where the concluding equality is given by $(b)$. We see that $x_t$ is feasible for (\ref{OptPrb}) and $\|x_t\|_1\leq \|w\|_1$ for all small positive $t$. Since $w$ is the unique optimal solution to (\ref{OptPrb}), we should have $x_t=w$, $t>0$, which would imply that $\mu_i=0$ for all $i$; but the latter is impossible by $(b)$. Thus, the optimal value in  (\ref{one}) is $<1$.
\par
We see that whenever $x$ is a vector with $s$ nonzero entries, equal to $\pm1$, there exists $y$ such that $(A^Ty)_i=x_i$ when $x_i\neq0$ and $|(A^Ty)_i|<1$ when $x_i=0$ (indeed, in the role of this vector one can take the $y$-component of an optimal solution to the problem (\ref{one}) coming from $w=x$), meaning that $\gamma_s(A)<1$, as claimed.
\par
b) Now assume that $\gamma_s(A)<1$, and let us prove that $A$ is $s$-good.
Thus, let $w$ be an $s$-sparse vector; we should prove that $w$ is the unique optimal solution to (\ref{OptPrb}). There is nothing to prove when $w=0$. Now let $w\neq0$, let $s'$ be the
number of nonzero entries of $w$, and $I$ be the set of indices of these entries. By {\bf C} we have $\gamma:=\gamma_{s'}(A)\leq\gamma_s(A)$, i.e., $\gamma<1$. Recalling the definition of $\gamma_s(\cdot)$, there exists $y\in\bbr^k$ such that $(A^Ty)_i=\sign(w_i)$ when $w_i\neq0$ and $|(A^Ty)_i|\leq\gamma$ when $w_i=0$. The function
$$
f(x)=\|x\|_1-y^T[Ax-Aw]=\sum_{i\in I}\left[|x_i|-\sign(w_i)(x_i-w_i)\right] +\sum_{i\not\in I}\left[|x_i|-\gamma_i x_i\right],
\,\gamma_i=(A^Ty)_i,\,i\not\in I,
$$
 coincides with the objective of (\ref{OptPrb}) on the feasible set of (\ref{OptPrb}). Since $|\gamma_i|\leq\gamma<1$, this function attains its  unconstrained minimum in $x$ at $x=w$. Combining these two observations, we see that $x=w$ is an optimal solution to (\ref{OptPrb}). To see that this optimal solution is unique, let $x'$ be another optimal solution to the problem. Then
$$
0=f(x')-f(w)=\sum_{i\in I}\underbrace{\left[|x^\prime_i|-\sign(w_i)(x^\prime_i-w_i)-|w_i|\right]}_{\geq0} +
\sum_{i\not\in I}\left[|x^\prime_i|-\gamma_i x^\prime_i\right];
$$
since $|\gamma_i|<1$ for $i\not\in I$, we conclude that $x^\prime_i=0$ for $i\not\in I$. This conclusion combines with the relation $Ax'=Aw$ to imply the required relation $x'=w$, due to the following immediate observation:
\begin{lemma}\label{kerA}
If $\gamma_s(A)<1$, then every $k\times s$ submatrix of $A$ has trivial kernel.
\end{lemma}
{\bf Proof.} Let $I$ be the set of column indices of a $k\times s$ submatrix of $A$. If  this submatrix has a nontrivial kernel there exists a nonzero $s$-sparse vector $z\in\bbr^n$
such that $Az=0$.  Let $I$ be the support set of $z$. By {\bf A}, there exists a vector $y\in\bbr^k$ such that $(A^Ty)_i=\sign(z_i)$ whenever $i\in I$, that is
 \[
 0=y^TAz=\sum_{i:z_i\neq0}(A^Ty)_iz_i=\|z\|_1,
 \]
 which is impossible. \epr
 \par
 (ii) Let $\gamma:=\gamma_s(A,\beta)<1$.  By definition it means that for every vector $z\in\bbr^n$ with $s$ nonzero entries, equal to $\pm1$, there exists $y$, $\|y\|_*\leq\beta$, such that $A^Ty$ coincides with $z$ on the support of $z$ and is such that $\|A^Ty-z\|_\infty\leq\gamma$. Given $z$, $y$ as above and setting $y'={1\over 1+\gamma}y$, we get $\|y'\|_*\leq{1\over 1+\gamma}\beta$ and \[
\|A^Ty'-z\|_\infty\leq\max\left[1-{1\over1+\gamma},{\gamma\over 1+\gamma}\right]={\gamma\over1+\gamma}.
 \]
 Thus, for every vector $z$ with $s$ nonzero entries, equal to $\pm1$, there exists $y'$ such that $\|y'\|_*\leq{1\over 1+\gamma}\beta$ and
$
\|A^Ty'-z\|_\infty\leq {\gamma\over 1+\gamma},$
 meaning that $\gamma:=\gamma_s(A,\beta)<1$ implies
\begin{equation}\label{meaningthat1}
 \widehat{\gamma}_s\left(A,{1\over1+\gamma}\beta\right)\leq{\gamma\over1+\gamma}<1/2.
\end{equation}
Now assume that $\widehat{\gamma}:=\widehat{\gamma}_s(A,\beta)<1/2$. For an $s$-element subset $I$ of the index set
$\{1,...,n\}$, let
$$
\Pi_I=\left\{u\in\bbr^n:\mbox{ exists }  y\in\bbr^k: \|y\|_*\leq\beta, \;(A^Ty)_i=u_i\,\mbox{ for } i\in I,\, |(A^Ty)_i|\leq\widehat{\gamma}\,\mbox{ for } i\in\bar{I}\right\},
$$
where $\bar{I}$ is the complement of $I$. It is immediately seen that $\Pi_I$ is a closed and convex set in $\bbr^n$. Let  $B$ be the centered at the origin
$\|\cdot\|_\infty$-ball of the radius $1-\widehat{\gamma}$ in $\bbr^n$: $B=\{u\in \bbr^n:\,\|u\|_\infty\le 1-\widehat{\gamma}\}$. We claim that $\Pi_I$ contains $B$. Using this fact we conclude that for every vector $z$ supported on $I$ with entries $z_i$, $i\in I$, equal to $\pm1$, there exists an $u\in\Pi_I$ such that $u_i=(1-\widehat{\gamma})z_i$, $i\in I$. Recalling the definition of $\Pi_I$, we conclude that there exists $y$ with $\|y\|_*\leq{(1-\widehat{\gamma})^{-1}}\beta$ such that $(A^Ty)_i={(1-\widehat{\gamma})^{-1}}u_i=z_i$ for $i\in I$ and $|(A^Ty)_i|\leq {(1-\widehat{\gamma})^{-1}}\widehat{\gamma}$ for $i\not\in I$. Thus, the validity of our claim would imply that
\begin{equation}\label{meaningthat2}
\widehat{\gamma}:=\widehat{\gamma}_s(A,\beta)<1/2\Rightarrow \gamma_s\left(A,{1\over1-\widehat{\gamma}}\beta\right)\leq{\widehat{\gamma}\over1-\widehat{\gamma}}<1.
\end{equation}
Let us prove our claim. Observe that by definition $\Pi_I$ is the direct product of its projection $Q$ on the plane $L_I=\{u\in\bbr^n:u_i=0,i\not\in I\}$ and the entire orthogonal complement $L_I^\perp=\{u\in\bbr^n:u_i=0,i\in I\}$ of $L_I$; since $\Pi_I$ is closed and convex, so is $Q$. Now, $L_I$ can be naturally identified with $\bbr^s$, and our claim is exactly the statement that the image $\bar{Q}\subset\bbr^s$ of $Q$ under this identification contains the centered at the origin $\|\cdot\|_\infty$ ball $B_s$, of the radius $1-\widehat{\gamma}$, in $\bbr^s$. Assume that it is not the case. Since $\bar{Q}$ is convex and $B_s\not\subset \bar{Q}$, there exists $v\in B_s\backslash\bar{Q}$, and therefore there exists a vector $e\in\bbr^s$,
$\|e\|_1=1$ such that $e^Tv>\max_{v'\in \bar{Q}}e^Tv'$ (recall that $Q$, and thus $\bar{Q}$, is both convex and closed). Now
let $z\in\bbr^n$ be the $s$-sparse vector supported on $I$ such that the entries of $z$ with indices $i\in I$ are the signs of the corresponding entries in $e$. By definition of $\widehat{\gamma}=\widehat{\gamma}_s(A,\beta)$, there exists $y\in \bbr^k$ such that $\|y\|_*\leq\beta$ and $\|A^Ty-z\|_\infty\leq\widehat{\gamma}$; recalling the definition of $\Pi_I$ and $\bar{Q}$, this means that $\bar{Q}$ contains a vector $\bar{v}$ with $|\bar{v}_j-\sign(e_j)|\leq\widehat{\gamma}$, $1\leq j\leq s$,  whence
$e^T\bar{v}\geq \|e\|_1-\widehat{\gamma}\|e\|_1=1-\widehat{\gamma}$. We now have
$$
1-\widehat{\gamma}\geq\|v\|_\infty\geq e^Tv>e^T\bar{v}\geq 1-\widehat{\gamma},
$$
where the first $\geq$ is due to $v\in B_s$, an $>$ is due to the origin of $e$. The resulting inequality is impossible, and thus our claim is true.
\par
We have proved the relations (\ref{meaningthat1}), (\ref{meaningthat2}) which are slightly weakened versions of (\ref{onehas1}.$a$-$b$). It remains to prove  that the inequalities $\leq$ in the conclusions of (\ref{meaningthat1}), (\ref{meaningthat2}) are in fact equalities. This is immediate: assume that under the premise of (\ref{onehas1}.$a$) we have \[
\widehat{\gamma}:=\widehat{\gamma}_s\left(A,{1\over 1+\gamma}\beta\right)<\gamma_+:={\gamma\over 1+\gamma}.
 \]
 When applying (\ref{meaningthat2})  with $\beta$ replaced with ${1\over 1+\gamma}\beta$, we get
\begin{equation}\label{tem1}
\gamma_s\left(A,{1\over 1-\widehat{\gamma}}\left[{1\over 1+\gamma}\beta\right]\right)\leq{\widehat{\gamma}\over 1-\widehat{\gamma}}<
{\gamma_+\over 1-\gamma_+}=\gamma.
\end{equation}
At the same time, ${1\over 1-\widehat{\gamma}}{1\over 1+\gamma}<{1\over 1-\gamma_+}{1\over 1+\gamma}=1$ due to $\widehat{\gamma}<\gamma_+$; since $\gamma_s(A,\cdot)$
is nonincreasing by {\bf B}, we see that
\[
\gamma_s\left(A,{1\over 1-\widehat{\gamma}}\left[{1\over 1+\gamma}\beta\right]\right)\geq \gamma_s(A,\beta),
\] and thus (\ref{tem1}) implies that
$
\gamma_s(A,\beta)<\gamma,
$
which contradicts the definition of $\gamma$. Thus, the concluding $\leq$ in (\ref{meaningthat1}) is in fact equality. By completely similar argument, so is the concluding $\leq$ in (\ref{meaningthat2}).  \epr
\end{document}